\documentclass{etds} 
\usepackage{amsfonts}
\usepackage{amsmath}
\usepackage{latexsym}
\usepackage{graphicx}
\usepackage{amssymb}
\usepackage{amsbsy}

\newtheorem{theorem}{Theorem}
\newtheorem{proposition}[theorem]{Proposition}
\newtheorem{lemma}[theorem]{Lemma}
\newtheorem{corollary}[theorem]{Corollary}

\def\endproof{\hfill{\vrule height4pt width6pt depth2pt}

\vspace{0.5em}

}
\def\btau{{\boldsymbol{\tau}}}
\def\ba{ {\boldsymbol{a}}}
\def\bv{ {\boldsymbol{v}}}
\def\bx{ {\boldsymbol{x}}}
\def\by{ {\boldsymbol{y}}}
\def\bu{ {\boldsymbol{u}}}
\def\bw{ {\boldsymbol{w}}}
\def\bp{ {\boldsymbol{p}}}
\def\bs{ {\boldsymbol{s}}}

\begin{document}\ETDS{1}{10}{28}{2008} 
\title{Branch points of substitutions and closing ordered Bratteli diagrams \footnote{This research was partially supported by an NSERC grant. }}

\runningheads{ R. Yassawi}{}
 \author{ Reem Yassawi}

\address{Dept. of Mathematics, Trent University \\
1600 West Bank Drive, Peterborough, Ontario, K9J 7B8, Canada.}

\email{ ryassawi@trentu.ca}

\recd{2009}





\begin{abstract} 
We study stationary ordered Bratteli diagrams and give necessary and
sufficient conditions for these orders to generate a continuous
Vershik map. We apply this to finding adic representations for one
sided substitution subshifts.  We give an algorithm to find the branch
points of a substitution, which have to be mapped to the minimal
elements of such an ordering. We find adic representations for
substitutions with one branch point, and also substitutions all of
whose branch points are fixed.

\end{abstract}

\section{Introduction}

While adic representations of primitive  two-sided
substitution subshifts have been studied in \cite{f} and  \cite{dhs}, and generalised to aperiodic two-sided substitution subshifts in 
\cite{bkm},  this article originated with the question of whether such
a representation existed for one-sided substitution subshifts. This
question is a special case of the general question of how far
two-sided concepts and results in \cite{hps} and \cite{gps} about
minimal homeomorphisms of a Cantor space can be generalised to minimal
continuous maps of a Cantor space.  Conjugacy of a pair of two-sided
subshifts does not imply conjugacy of the corresponding one-sided
subshifts (see Example \ref{notonesidedconjugate} for an example where
the subshifts are generated by substitutions).

An example of a markedly different situation in the one-sided case is
the generalisation of Theorem 4.6 in \cite{hps}, which tells us that any
minimal homeomorphism $(X,T)$ of a Cantor space has an adic
representation. The proof of this result requires a sequence of
refining, generating clopen partitions $({\cal P}_{n})$, which naturally define an ordered Bratteli diagram ${\cal B}$. If the 
 intersection of the  bases of ${\cal P}_{n}$ is a unique point, then the ordering is {\em proper}, so that the Vershik-adic map 
$V_{\cal B}:X_{\cal B}\rightarrow X_{\cal B}$ representing $T$ is
a homeomorphism. For a continuous minimal map $T:X\rightarrow X$
of a Cantor space, first there is the question of finding the
non-invertible points for $T$- we call these {\em branch points};
second we need a sequence of refining, generating clopen partitions
the intersection of whose bases is the branch points.  In addition
though the ordered Bratteli diagram generated by this sequence of
partitions must have the correct number of maximal paths, and
must also  support a continuous adic map, so, in the case where there are $n$ branch points, there is the
question of determining when an ordered Bratteli diagram with $n$
minimal paths and $N$ maximal paths has a continuous Vershik map
defined on it - we call such ordered diagrams {\em closing}.

For (aperiodic and  primitive) substitution subshifts $(X_{\btau},\sigma)$ -see Section \ref{preliminaries} for definitions -  the natural sequence of generating partitions $({\cal P}_{n})$ to consider is that where ${\cal P}_{n}$ is  
generated by the ${\btau}^{n}$-words. The resulting ordered Bratteli diagram ${\cal B_{\btau}}$ has as many minimal elements as there are right fixed $\btau$-fixed points $ x_{0}x_{1}\ldots$
in $X_{\btau}^{\mathbb N}$ and as many maximal elements as there are left fixed points 
$\ldots x_{-2}x_{-1}$ in $X_{\btau}^{\mathbb N}$.
The proof in \cite{dhs} of the statement that every two sided
(aperiodic, primitive) substitution subshift $(X_{\btau}^{\mathbb Z},
\sigma)$ is topologically conjugate to some adic system $(X_{\cal
B},V_{\cal B})$ consists then of 
finding a substitution $\btau'$ such that $(X_{\btau}^{\mathbb Z},
\sigma) \cong (X_{\btau'}^{\mathbb Z}, \sigma)$ and where $\btau'$ is
{\em proper} -ie has exactly one of each of these left and right $\btau$-fixed points.

Suppose $(X_{\btau}^{\mathbb N},\sigma)$ has $n$ branch points $\bx_{1},\ldots \bx_{n}$ where $\bx_{i}$ is a $k_{i}$-branch point, (so
that $\bx_{i}$ has $k_{i}$ $\sigma$-preimages).
Finding an adic representation for $(X_{\btau}^{\mathbb N},\sigma)$, the one-sided subshift defined by $\btau$, requires
finding $\btau'$ such that

\begin{equation}
\label{adiccondition0}
(X_{\btau'}^{\mathbb N},\sigma) \stackrel{\Phi}{\cong} (X_{\btau}^{\mathbb N},\sigma),
\end{equation}

\begin{equation}
\label{adiccondition1}
\btau' \mbox{ has }n \mbox{ right fixed points } 
\{{\bf v}_{i}\}_{i=1}^{n} \mbox{  with }
\Phi({\bx}_{i})={\bf v}_{i} \mbox{ for } 1\leq i \leq n,
\end{equation}

\begin{equation}
\label{adiccondition2}
\mbox{ $\btau'$ has $\sum_{i=1}^{n}k_{i}$ left fixed points 
$\{{\bu}_{j}^{i}:\,\,\,1\leq j\leq k_{i}, \,\,\, 1\leq  i\leq n 
\}$,} 
\end{equation}
and
\begin{equation}
\label{adiccondition3}
\mbox{ $u_{j}^{i}v_{i'}\, \in {\cal L}_{\btau'}$ if and only if $i=i'$}.
\end{equation}

Here ${\cal L}_{\btau'}$ is the language of $\btau'$, and the  right  (left) fixed point $\bu$ is generated by the letter $u$ whose substitution word has $u$ as a prefix (suffix). Condition \ref{adiccondition1} is a
requirement that branch points of $\btau$ are mapped to fixed points
of $\btau'$, so that the  sequence of generating
partitions defined by $\btau'$ have bases whose intersection is $\{{\bf v}_{1},\ldots ,{\bf v}_{n}\}$, and also 
span the topology of
$X_{\btau'}^{\mathbb N}$ (Proposition \ref{clopenpartition}).
 Condition \ref{adiccondition2}   ensures that an adic representation of $\btau'$ has the right number of maximal elements.
 Condition \ref{adiccondition3} is the key to ensuring that the ordered diagram ${\cal B}_{\btau'}$ is closing.
 This is explained in Proposition
\ref{closingproposition}.

In Section \ref{preliminaries} we lay out definitions and state standard results on substitutions that will be used.
We assume throughout that $\btau$ is aperiodic, recognizable and all
of its powers are injective on letters.  Apart from Example
\ref{chacon}, we assume $\btau$ is primitive. These are standard
assumptions so that Proposition \ref{clopenpartition} applies.
 An algorithm which finds branch points of a substitution is described
 completely in Section \ref{suffix}. We remark that branch points are
 the one-sided counterpart to the {\em left (or right) asymptotic}
 points in a two-sided subshift, which have been studied, for example,
 in \cite{lzpre} and \cite{BDH} - two bi-infinite points are right
 asymptotic if they agree from some index $N$ onwards - for us the
 right infinite ray of the $N$-th shift of these points is a branch point.  
We characterise non-$\btau$-fixed left branch points $\by$ of $\btau$ 
 as points which satisfy an equation of the form ${\bf
w}\btau(\by)=\by $ for some non-empty finite word ${\bf
w}$    (Lemma \ref{nonfixedlemma1a}). If $\bu$ is a $\btau$-fixed point, Proposition \ref{inverseimagesofu} characterises when $\bu$ is a branch point.

In Section \ref{Bratteli}  we define stationary Bratteli diagrams and prove Proposition
\ref{closingproposition}, which characterises in terms of ${\cal L}_{\btau}$ when an adic
map on a stationary diagram $B$ with an ordering generated by $\btau$
is continuous. This result has a generalisation to arbitary finite
rank ordered diagrams and this is the basis for most of the results  in  \cite{BKMY}.

In Section \ref{quasi_invertible} we deal with the special case where the substitution
subshift $(X_{\btau}^{\mathbb N},\sigma)$ has a unique branch point, and
call these subshifts {\em quasi-invertible}. Subshifts with this
property have been studied: Sturmian sequences, and more generally,
 Arnoux-Rauzy sequences (\cite{AR}) generate quasi-invertible
subshifts - the condition (*) in \cite{AR} ensures this, and in many
cases these sequences are fixed points of substitutions (see \cite
{cmps} for a characterization of when Sturmian subshifts are
substitution subshifts; if the `S-adic' expansion of a sequence in
\cite{AR} is periodic, then it is the fixed point of a substitution,
for example the tribonacci substitution). All these subshifts have
nice geometric representations as toral rotations, or interval
exchange maps.  Many (but not all) Pisot substitutions that have been
studied in the literature are quasi-invertible. Any substitution with
a rational Perron eigenvalue has the two-sided subshift 
$(X_{\btau}^{\mathbb Z},\sigma)$ orbit equivalent to $(X_{\btau'}^{\mathbb Z},\sigma)$
where $\btau'$ is
quasi-invertible  (see Example \ref{orbitequivalence}).
It would be interesting to study whether quasi-invertibility imposes
any spectral constraints. While not all have discrete spectrum (see
Example \ref{chacon}: the Chacon substitution is quasi-invertible)
spectral multiplicity may be bounded.
If $\btau$ is quasi-invertible, then Condition (\ref{adiccondition3}) is immediately 
satisfied, and that a substitution satisfying Conditions
(\ref{adiccondition0}), (\ref{adiccondition1}) and (\ref{adiccondition2}) exists 
is the statement of Theorems  \ref{quasiinvertibleleftproperareadic} and 
\ref{quasiinvertiblearealmostadic}. We remark  that a $\btau'$ satisfying Conditions 
(\ref{adiccondition0}) and (\ref{adiccondition1}) for 
constant length substitutions on a two-letter alphabet is found in Section 7 of
\cite{CK}.

 If $\btau$ has several branch points, all of which are
fixed, then we prove in Theorem \ref{generalfixedbranch} that the
right $\btau'$ can be found. 
For general $\btau$ we are not able to find
the right substitutions that satisfy Condition \ref{adiccondition1} and we pose the question of whether given a $\btau$, $(X_{\btau}^{\mathbb N},\sigma)\cong (X_{\btau'}^{\mathbb N},\sigma)$ where all of $\btau$'s branch points are mapped to $\btau'$-fixed points.

The author thanks Ian Putnam and Fabien Durand for helpful discussions.

\section{Preliminaries \label{preliminaries}}
\subsection{Notation} 

 Let ${\cal A}$ be a finite alphabet.  A {\em word} $\bw$ (denoted using
 boldface) from ${\cal A}$ is a finite concatenation of elements from
 ${\cal A}$.  Let ${\cal A}^{+}$ denote the set of words (including
 the empty word) from ${\cal A}$.  If ${\bw}=w_{1}\ldots w_{m}$ and
 ${\bw^{*}}=w_{1}^{*}\ldots w_{n}^{*}$, then ${\bw\bw^{*}}:= w_{1}
 \ldots w_{m}w_{1}^{*}\ldots w_{n}^{*}$. Define $\bw_{[i,j]}=
 w_{i}w_{i+1}\ldots w_{j}$.  If $\bw= w_{1}\ldots w_{m}$, the {\em
 length} of $\bw$, written $|{\bw}|$, is $m$. Given words $\bw$,
 ${\bw^{*}}$, say that $\bw$ {\em occurs} in ${\bw^{*}}$, or $\bw$
 {\em is a subword of} ${\bw^{*}}$, if $\bw= {\bw^{*}}_{[i,
 i+|{\bw}|-1]}$ for some $i$.  If $|\bw|=m$, define the {\em k-shift}
 $\sigma^{k}(\bw)=\bw_{[k+1, m]} $. ${\bw^{*}}$ is a {\em prefix} of
 $\bw$ if $\bw_{[1,|{\bw^{*}}|]}= {\bw}^{*}$, similarly ${\bw^{*}}$ is
 a {\em suffix} of $\bw$ if $\bw_{[|\bw|-|{\bw^{*}}|+1,|\bw|]}=
 \bw^{*}$. When we speak of `the' prefix (suffix) of a word $\bw$, we
 will mean the first (last) letter of $\bw$.

 Let ${\mathbb N}:=\{0,1, \ldots\}$. If ${\mathbb M} = {\mathbb Z}$ or
 ${\mathbb N}$, then the space of all ${\mathbb M}$-indexed sequences
 from ${\cal A}$ is written as ${\cal A}^{\mathbb M}$, and a
 configuration $\bx\,\in {\cal A}^{\mathbb M}$ is written $\bx=
 (x_{m})_{m\, \in {\mathbb M}}$. Let ${\cal A}$ be endowed with the
 discrete topology and ${\cal A}^{\mathbb M}$ with the product
 topology (which is metrizable);
if $b\,\in {\cal A}$
 and $j\, \in \mathbb M$, the clopen sets $[b]_{j}:=\{\bx:
\,\,x_{j}=b\}$ form a countable basis for the
topology on ${\cal A}^{\mathbb M}$ -if $j=0$; we write $[b]$ instead of  $[b]_{0}$.   The (left) {\em shift map} $\sigma: {\cal
A}^{\mathbb M} \rightarrow {\cal A}^{\mathbb M}$ is the map defined as
$(\sigma(\bx))_{m} = x_{m+1}$.  $(X, \sigma)$ is a {\em subshift of $({\cal
A}^{\mathbb M},\sigma)$} if $X$ is a closed $\sigma$-invariant subset
of ${\cal A}^{\mathbb M}$. 


\subsection{Substitutions}

A {\em substitution} is a map $\btau: {\cal A} \rightarrow {\cal
A}^{+}$.
We extend $\btau$ to a map $\btau: {\cal
A}^{+}\rightarrow {\cal A}^{+}$ by concatenation: if ${\ba}=
a_{1}\ldots a_{k}$, then $\btau({\bf a}):= \btau(a_{1})
\ldots \btau(a_{k})$. In this way iteration $\btau^{n}$ is well defined. 
The substitution $\btau$ is extended to a map $\btau: {\cal A}^{\mathbb
N} \rightarrow {\cal A}^{\mathbb N} $ defined by $ \btau (\bx) := 
\btau (x_{0}) \,\btau (x_{1})\ldots $, and also  $\btau: {\cal A}^{\mathbb
Z} \rightarrow {\cal A}^{\mathbb Z} $ defined by $ \btau (\bx) :=
\ldots \btau (x_{-1}) \cdot \btau (x_{0}) \,\btau (x_{1})\ldots $. 
 We say $\btau$ is {\em left (right) proper} if there exists $l$ in ${\cal A}$
such that $l$ is the  prefix (suffix) of $\btau (a)$ for each $a\,\in {\cal A}$.
 If $\btau$  is both left and right
proper, it is called {\em proper}. We say ${\btau}$ {\em suffix-permutative} if the set of suffixes 
 of $\{\btau(a)\}_{a\in {\cal A}}$ is a permutation of ${\cal
A}$. Also $\btau$ is {\em injective} if $\btau(a)\neq \btau(b)$ whenever $a\neq b$.
If $|{\cal A}|=d$, then let 
$M=(m_{i,j})_{1\leq i,j\leq d}$ be the {\em incidence matrix} for ${\btau}$, 
where
 $m_{i,j}$ is the number of occurences of $j$ in $\btau(i)$.
We say ${\btau}$ is {\em primitive} if there exists a
positive integer $k$ such that for any $a\, \in {\cal A}$, all letters
of ${\cal A}$ appear in $\btau^{k}(a)$, ie some power of the incidence matrix is positive. Let  ${\cal L}_{\btau}$ be the language generated by the words
$\{\btau^{n}(a): n\, \in {\mathbb N}^{+},\,\,a
\, \in {\cal A}\}$.
  A (right infinite) {\em fixed point} of $\btau$ is a sequence ${\bu
  } \, \in {\cal A}^{\mathbb N}$ such that $\btau({\bu}) = {\bu}$. If
  ${\bu}\, \in {\cal A}^{\mathbb N}$ is a fixed point for $\btau$,
  then $\btau(u _{0})$ starts with $u_{0}$.  Conversely, if there
  exists a letter $l$ such that $l$ is a prefix of $\btau(l)$, then
  ${\bu}:=\lim_{n\rightarrow \infty} \btau^{ n}(l)$ is the unique
  fixed point satisfying $u_{0}=l$.  Similarly one can define a {\em left
  infinite} fixed point using a letter $u_{0}$ such that $\btau(u_{0})$
  ends with $u_{0}$.  Using the pigeonhole principle, there exists
  $n\geq 1$ such that $\btau^{ n}$ has at least one fixed point. If
  $\btau$ is primitive, then any fixed point $\bu$ generates the same
  subshift $X_{\btau}^{\mathbb N}:= \overline{\{\sigma^{n}({\bu})\}}$.
  Henceforth we will be working mainly with one sided substitution
  subshifts, and unless otherwise indicated, will assume this, and
  write $(X_{\btau},\sigma)$ instead of $(X_{\btau}^{\mathbb
    N},\sigma)$.

In what follows we will almost always assume that our substitution
$\btau$ is aperiodic and primitive, the exception being in Example
\ref{chacon}.  There $\btau$ has a fixed point $\bu$ which generates
${\cal L}_{\btau}$ (we call this a {\em generating} fixed point) and
the resulting subshift $(X_{\btau},\sigma)$ is still minimal (ie for
each $x$ in $X_{\btau}$, $\{ \sigma^{n}(x): n \geq 0\}$ is dense in
$X_{\btau}$).  Such substitutions enjoy many established properties
that primitive substitutions do, and minor modifications of proofs of
facts for primitive substitutions can be made to ensure that the same
results hold for these substitutions -see \cite{bkm}.

{\em Unilateral recognizability} was a condition introduced in
\cite{Ho} and \cite{Qu} which ensured that one could uniquely `de-substitute' sequences in $X_{\btau}$, and  find clopen
generating partitions for $(X_{\btau},\sigma)$. Specifically, the set of
(1)-cuttings of $\bu$ are
\[E:=\{0\}\cup\{ |\btau (\bu_{[0,p-1]})|: p>0\}\,.\]
We then say that $\btau$ is {\em (unilaterally) recognizable} if there
exists some $L$ so that if $\bu_{[i,i+L-1]}= \bu_{[j,j+L-1]}$ and
$i\, \in E$, then $j\, \in E$. A word $\bw$ occurring at ${\bf
u}_{i}$ and $\bu_{j}$ has the {\em same 1-cutting} at $i$ and $j$,
with $i<j$, if $E\cap \{i,\ldots i+|\bw|\} = j-i + E\cap \{ j ,\ldots j+|\bw|\}.$

The necessary and sufficient
conditions that exist in Theorem 3.1, \cite{Mo1} for a primitive substitution to be unilaterally recognizable 
can be modified (see \cite{bkm}) to work for  substitutions:

\begin{theorem}
\label{mosserecognizability}
Let $\btau$ be minimal  with generating  fixed point $\bu$. $\btau$ is not
recognizable if and only if for each $L$, there exists a word $\bw$ of
length L, and two elements $a, \, b \,\in {\cal A}$ such that 
\begin{enumerate}
\item $\btau(b)$ is a proper suffix of $\btau(a)$, and 
\item The words $\btau(a) \bw$ and $\btau (b) \bw$ appear in 
$\bu$ with the same 1-cutting of $\bw$.
\end{enumerate} \endproof
 \end{theorem}

 In particular,
suffix permutative substitutions are recognizable. It is
straightforward to show that if $\btau$ is recognizable, and injective,  then 
$\btau^{n}$ is recognizable for each $n$.

In Lemmas 2 and 3, \cite{Ho}, the following is shown for primitive
recognizable substitutions (a similar proof works if $\btau$ is minimal).   Note that $\btau (X_{\btau}) \subset X_{\btau}$.

\begin{proposition}\label{hosttopologyprop}

 Let $\btau$ be a primitive,  recognizable  substitution with 
generating fixed point $\bu$. 

\begin{enumerate}
\item
$\by\, \in \btau (X_{\btau})$ if and only if $\sigma^{n_i} \bu
\rightarrow \by$ where $n_{i}\, \in E$ for all large $i$.
\item
$\btau ({X}_{\btau})$ is clopen in $X_{\btau}$.
\item
$\sigma^{p}(\btau (x)) \ \in \btau(X_{\btau})$ if and only if $p=|\btau(\bx_{[0,r]})|$ 
for some $r$ so that $\sigma^{p}(\btau(x)) = 
\btau(\sigma^{r+1}(\bx))$.
\end{enumerate}
If $\btau$ is also injective, then
\begin{enumerate}
\setcounter{enumi}{3}
\item
$\btau:X_{\btau}\rightarrow \btau(X_{\btau} )$ is a homeomorphism. \endproof
\end{enumerate}
\end{proposition}

We remark here that the previous proposition also implies that for each $n\geq 1$, $\btau^{n}:X_{\btau}\rightarrow {\btau}^{n}(X_{\btau})$ is also a homeomorphism, though as a substitution $\btau^{n}$ is not necessarily injective. 

Henceforth we assume that our substitution $\btau$ is 
recognizable  and all powers of $\btau$ are injective, though the latter assumption is
not required for the section on branch points.

\begin{corollary}
\label{almostuniquerep}

Every $\by\,\in \btau (X_{\btau})$ can be written in a unique way as ${\bf
y}= \btau (\bx)$ with $\bx\,\in X_{\btau}$. Every $\by \, \in X_{\btau}
\backslash \btau (X_{\btau})$ can be written as $\by= \sigma^{k} \btau
(\bx)$ where $\sigma(\bx)$ is unique, and $0<k<|\btau (x_{0})|$.

\end{corollary}

{\bf Proof:} The first assertion follows immediately from Proposition
\ref{hosttopologyprop}. If  $\by \,\in X_{\btau}\backslash \btau (X_{\btau})$ then by 
Part 1 of Proposition \ref{hosttopologyprop}, there is a sequence
$\{n_{i}\}$ of integers not in E, with $\sigma^{n_i}(\bu) \rightarrow
\by$.  Now each $n_{i}= k_{i} + r_{i}$ where $k_{i} \, \in E$, 
and the positive 
$r_{i}$ is strictly less than the successor of $k_{i}$ in $E$. Since
$E$ is almost periodic, by dropping to a subsequence if necessary, we
assume that there is some fixed positive $r$ such that each 
$n_{i}= k_{i} +r$. Now we can assume that $\sigma^{k_i}(\bu)$ converges to
some $\bx \, \in \btau(X_{\btau})$. So $\by = \sigma^{r}\bx$.

Now suppose that $\by = \sigma^{r} \bx = \sigma^{r'}\bx'$, where  
if $\bx = x_{0}x_{1}\ldots$ then $r< |\btau(x_{0})|$, and similarly for 
$r'$. Let $l= |\btau (x_{0})|-r$, and define $l'$ similarly.

If $l=l'$, then $\sigma^{r+l} \btau(\bx)= \btau (\sigma\bx)\, \in
\btau (X_{\btau})$ and similarly $\sigma^{r'+l} \btau (x') = \btau (\sigma
\bx') \, \in \btau (X_{\btau})$, so $\btau (\sigma \bx) = \btau (\sigma
\bx')$. Thus $\sigma \bx = \sigma \bx'$.

Finally we will show that $l\neq l'$ leads to a contradiction. Suppose
$l<l'$. Then $\sigma^{r+l} \btau(\bx)= \btau (\sigma\bx)\, \in \btau
(X_{\btau})$, and so $ \sigma^{r'+l} \btau(\bx') = \sigma^{r+l} \btau(\bx)\in 
\btau (X_{\btau})$. Thus by part 3 of Proposition \ref{hosttopologyprop}, $r'+l = 
|\btau (\bx'_{[0,j]})|$ for some $j$. This contradicts  
$r'+l< r'+l' = |\btau (x'_{0})|$.
\endproof

\section{Branch points \label{suffix}}

If $(Y,T)$ is a system where $\by$ is not $T$-invertible, we
call $\by$ a {\em branch point}; if $\by$ has $M$ pre-images under $T$, we
say that ${\bf y}$ is an $M$-branch point. Let us call systems $(Y,T)$
with only one branch point {\em quasi-invertible}.  If
$(X_{\btau},\sigma)$ is quasi-invertible, and the branch point is an
$M$-branch point, we'll call $\btau$ $M$-{\em quasi-invertible}. More generally, if a substitution $\btau$ has $n$ $\btau$-right fixed points
$\{\bx_{1},\ldots \bx_{n}\}$ where $\bx_{i}$ is a $k_{i}$-branch point
and $N:=\sum_{i=1}^{n}k_{i}$, call $\btau$ an {\em
$(N,n)$-substitution}.
In this section we describe the substitution structure of branch points $\bx\in X_{\btau}$, and describe an algorithm which finds all such points. Say that $a_{1}$ is {\em part of a
suffix cycle} $a_{1}, a_{2},\ldots a_{m}$ if
$\btau(a_{i})$ has $a_{i+1}$ as a suffix for $1\leq i<m$, and $\btau (a_{m})$ has $a_{1}$ as a suffix.

\begin{proposition}
\label{inverseimagesofu}
Suppose that $\btau$ has a generating fixed point $\bu =
u_{0}u_{1}u_{2}\ldots$.  Then $a\bu\, \in X_{\btau}$ if and only if
$a$ is part of a suffix cycle $a_{1},a_{2},\ldots a_{m}$, and
$a_{i}u_{0}\, \in {\cal L}_{\btau}$ for some $i\, \in \{1,\ldots m\}$.
\end{proposition}

{\bf Proof:} If for some $i$, $a_{i}u_{0}\, \in {\cal L}_{\btau}$,
then since $\btau$ is primitive, there exists an $n$ such that
$a_{i}u_{0}\,\in {\btau}^{n}(u_{0})$, so $a_{i+1}\,\btau(u_{0})\, \in
\btau^{n+1}(u_{0})$, and $a_{i}\btau^{nm}(u_{0})\,\in
     {\btau}^{nm}(u_{0})$ for each $n\geq 1$. It follows that
     $a_{i}\bu\, \in X_{\btau}$, and that this is true for each
     $a_{j}\, \in \{a_{1},\ldots a_{m}\}.$

Conversely, suppose that $a\,\bu \, \in X_{\btau}$. We claim first
that $a$ is the suffix of some substitution word. If $a=\btau(b)$ for
some $b$, we are done. If not, by Corollary \ref{almostuniquerep},
$a\bu={\bw}\btau(\sigma \bx)$ where ${\bw}$ is the suffix of some
substitution word $\btau(x_{0})$, and $\sigma \bx$ is unique. Thus
$\btau(\bu)= \bu = \sigma ({\bw})\btau(\sigma \bx)$, and by Corollary
\ref{almostuniquerep} again, $\sigma (\bw) = \btau (u_{0}\ldots
u_{r})$ for some $r\geq 0$, so that $a$ is the suffix of a
substitution word.

Since $(X_{\btau}, \sigma)$ is minimal, $\sigma$ is onto, which means
that for each $n$ there are words $\bw_{n}$ of length $n$ such that
$\bw_{n}a\bu \, \in X_{\btau}$. Given $N$, choose $n$ so that $n>
|\btau^{N}(\alpha)|$ for all $\alpha$.  Thus $\bx = \bw_{n}a \bu =
\bw'\btau^{N}(\sigma (\by))$ where $\bw'$ is the suffix of
$\btau^{N}(y_{0})$, and Corollary \ref{almostuniquerep} implies that
for some $p>0$, $\sigma^{p+1} (\by) = \bu$ and $\bw_{n}a = \bw'
\btau^{N}(\by_{[1,p]})$, so that $a$ is the suffix of a $\btau^{N}$
word. Taking a subsequence, we conclude that there is some $b$ such
that $a$ is the suffix of $\btau^{N}(b)$ for infinitely many
$N$'s. Since the sequence of suffixes of $\btau^{N}(b)$ is eventually
periodic, this implies that $a$ is part of a suffix cycle.  \endproof

It follows that if $\btau$ is left proper and suffix permutative, then
$\bu$ has $|{\cal A}|$ $\sigma$-preimages, and if $\btau $ is left and
right proper, then $\bu$ has one $\sigma$-preimage.

Let ${\cal W}_{n} :=\{ \bw: \bw=\sigma^{k}(\btau^{n}(a)):a \in {\cal A},\,\,
 0\leq  k<|\btau^{n}(a)|\}$.
If $\bw \, \in {\cal W}_{n}$ and 
$\bw=\sigma^{k}(\btau^{n}(a))$ with $k\geq 1$
, let
$r(\bw):=|\{\alpha: \alpha\bw \, \in {\cal W}_{n}\}|$, and let
 $\sigma^{-1}\bw= \{\alpha\bw: \alpha\bw \, \in {\cal W}_{n}\}$.

\begin{proposition}
\label{sigmainjectivity}
   Suppose that $\by$ is not $\btau$-fixed, and  $\by\, \in
 [\bw]$ with $\bw \in {\cal W}_{n}$. Then $\by$ is a branch point 
only if $r(\bw)>1$.
Conversely, suppose that there exists a sequence of words $\bw_{n}\,
\in {\cal W}_{n}$ with $|\bw_{n}|\rightarrow \infty$, with $\cap_{n}[\bw_{n}] = \by$.
If  $r({\bw}_{n})>1$ for each $n$,
 then $\by$ is a branch point.

\end{proposition} 

{\bf Proof:} Since $\by$ is not fixed, then for some $N$, $\by\not\in
\btau^{n}(X_{\btau})$ for each $n\geq N$, so that $\by=
     {\bw}_{n}\btau(x)$ with $\bw_{n} \in {\cal W}_{n}$. Then $\by$ is
     a branch point only if $r({\bf w_{n}})>1$ for each $n\geq N$.

If $\by=\cap_{n}[{\bw}_{n}]$ where $r({\bw}_{n})>1$ for each $n$, then
     passing down to a subsequence if
     necessary, we can assume that there are distinct $\alpha$ and
     $\beta$ so that $\alpha \bw_{n}$ and $\beta \bw_{n}$ are elements
     of ${\cal W}_{n}$, hence ${\cal L}_{\btau}$. Thus $\alpha \by$
     and $\beta \by$ are both elements of $X_{\btau}$, so that $\by$
     is a branch point.  \endproof

 Proposition \ref{sigmainjectivity} suggests an algorithm for
 generating all non-$\sigma$ invertible points, and checking whether
 they are the $\btau$-fixed point. We develop this algorithm in what
 follows, as our aim is to find the branch points of  $\btau$.
If $\btau$ is defined
 on ${\cal A}=\{a_{1},\ldots, a_{n}\}$, pick ${A}_{1}\subset {\cal
   A}$, where ${\bf s}_{1}$, the maximal proper common suffix of
 $\{\btau (a)\}_{a\, \in { A_{1}}}$, is non-empty. Thus for each $a$
 in ${A_{1}}$, $\btau (a) = {\bw}_{a}x_{a}{\bf s}_{1}$, and the set
 $A_{2}:=\{x_{a}:a\in A_{1}\}$ has at least two distinct
 elements. Call the set $A_{2}$ the set of ${\bf s}_{1}$-{\em
   predecessors} in $\{\btau(a):a\, \in A_{1}\}$. If $\{\btau(a): a\,
 \in A_{2}\}$ has ${\bf s}_{2}$ as a maximal proper common suffix, let
 $A_{3}$ be the set of ${\bf s}_{2}$-predecessors in $\{\btau(a): a \,
 \in A_{2}\}$.  Inductively, if $A_{k}$ is the set of ${\bf
   s}_{k-1}$-predecessors in $\{\btau(a): a\,\in A_{k-1}\}$, let ${\bf
   s}_{k}$ be the maximal proper common suffix of $\{\btau(a): a\, \in
 A_{k}\}$, and $A_{k+1}$ the set of ${\bf s}_{k}$-predecessors in
 $\{\btau(a): a\in A_{k}\}$.  This process may end after a finite
 number of steps, ie $A_{k}=\emptyset $ for some $k$ - in which case
 we discard this $A_{1}$ and start over with another subset of ${\cal
   A}$.  The sequence $(A_{k})$ is eventually periodic. For, the
 cardinality $|A_{k}|$ of the sets $A_{k}$ cannot increase. If
 $|A_{k}|$ stays constant as $k$ increases, then for some $k$ and $m$,
 $A_{k}=A_{k+m}$, and so ${\bf s}_{k}= {\bf s}_{k+m}$, which in turn
 implies that $A_{k+1}= A_{k+m+1}$, and by induction that $A_{k}$ is
 eventually periodic. If $|A_{k}|$ decreases, then repeat this
 argument. Thus $({\bf s}_{k})_{k\,\in {\mathbb N}}$ is also
 eventually periodic. Note that there can be more than one choice for
 ${ A}_{1}$ which leads to a non-trivial sequence of ${\bf s}_{k}'s$
 and $A_{k}$'s, and that there are only finitely many sequences $({\bf
   s}_{k})_{k\,\in {\mathbb N}}$ that can be generated, as ${A}_{1}$
 ranges over all subsets of ${\cal A}$.

We illustrate this terminology with the
following examples.
\begin{enumerate}

\item 
If $\btau(a)= aab$, and $\btau(b)=abb$, if $A_{1}={\cal A}$, then
$A_{n}={\cal A}$ and ${\bf s}_{n}=b$ for all $k\geq 2$. (In fact if
$\btau$ is any  substitution on a two letter alphabet whose
substitution words have a proper maximal non-empty suffix, then it is
always the case that $({\bf s}_{k})_{k\geq 1}$ is constant.)

\item
If $\btau(a)=acb$, $\btau(b)=aba$ and $\btau(c)=aaa$, then the only
possible choice for ${A}_{1}$ is ${ A}_{1}=\{b,c\}$. Here ${\bf
s}_{1}= a$, ${ A}_{2}=\{a,b\}$, and ${\bf s}_{2} = \emptyset$ as
$\btau (a)$ and $\btau (b)$ have no common suffix. Thus 
all  choices of 
${A}_{1}$  lead to sequences $({\bf s}_{n})$ which are eventually empty.

\item
If $\btau(a)=acb$, $\btau(b)=aba$ and $\btau(c)=aca$, letting ${A}_{1}=\{b,c\}$, 
then ${ A}_{k}= {A}_{1}$ and ${\bf s}_{k}={\bf s}_{2}$ for all $k\geq 2$.
\end{enumerate}

\begin{lemma}
\label{nonfixedlemma1a}
Let ${\btau}$ be defined on ${\cal A}=\{a_{1},\ldots ,a_{n}\}$, and,
given $A_{1}\subset {\cal A}$ , $({\bf s}_{m})_{m=1}^{\infty}$ and
$(A_{m})_{m=1}^{\infty}$ as defined and not eventually empty. Suppose
that $ {\bf s}_{k+p} = {\bf s}_{k}$ for $k\geq M$. Then
\begin{enumerate}
\item
For $1\leq n\leq M+p-1$, the $\btau^{n}$-substitution words have a
maximal proper common suffix ${\bf
S}_{n}=\left(\prod_{k=0}^{n-1}{\btau}^{k}({\bf s}_{n-k})\right)$, and
for $n=M+kp+j$, where $k\geq 1$ and $0\leq j\leq p-1$, the
$\btau^{n}$-substitution words have a maximal proper common suffix
${\bf S}_{n}={\bf P}_{n}\,{\bf I}_{n}$, where \begin{equation}
\label{yuk}{\bf P}_{n}=\left(\prod_{m=0}^{j}\btau^{m}({\bf s}_{M+j-m})\right)
\prod_{l=0}^{k-1}\left(\prod_{i=0}^{p-1}
\btau^{n-(M+(k-1-l)p+p-1-i)}({\bf s}_{M+(p-1-i)})\right)
\end{equation}
and 
\begin{equation} \label{yuk*}{\bf I}_{n}   =\left(        \prod_{m=1}^{M-1}{\btau}^{n-(M-m)}({\bf s}_{M-m})\                                                                                            \right) 
.
\end{equation}
\item
If $({\bf s}_{n})$ is eventually the  $k$-periodic sequence $({\bf w}_{i})_{i=1}^{k}$, then it defines $k$ branch points, with each branch point ${\by}$ satisfying an equation of the form 
\[({\bf w}_{i}\btau({\bf w}_{i+1})\btau^{2}({\bf w}_{i+2})\ldots 
\btau^{k-1}({\bf w}_{i-1}))\btau({\by}) ={\by}\,.\]
Also any branch point ${\by}$ which is non-$\btau$-fixed arises in this fashion.
\item
The substitution $\btau$ is quasi-invertible, with branch point $\by$
non $\btau$-fixed, if and only if all $\btau$-fixed points are
invertible, and there is an ${\bf s}\in {\cal A}^{+}$ such that any
$A_{1}\subset {\cal A}$ generates a sequence $({\bf s}_{n})$ which is
either eventually empty, or eventually the constant sequence $({\bf
s})$.

\item
 The substitution $\btau$ is quasi-invertible, with branch point ${\bf u}$, a $\btau$-fixed point, if and only  ${\bf u} $ is the only $\btau$-fixed branch point, and 
for all  subsets ${\cal A}_{1} \subset {\cal A}$ that generate a sequence $({\bf s}_{n})$ which is not eventually empty , there is
some $N$ such that for $n>N$, there is a sequence $(k_{n})_{n\geq N}$
of natural numbers such that the maximal common proper suffix of the
$\{\btau^{n}(a):a\,\in A_{1}\}$ have $\btau^{k_{n}}(u_{0})$ as a prefix,
where $k_{n}\rightarrow \infty$.
\end{enumerate}
\end{lemma}

{\bf Proof}: 
\begin{enumerate}
\item
The proof of Equations \ref{yuk} and \ref{yuk*} follow from the definitions of $({\bf s}_{k})$, and induction.
\item
Note that if $\lim_{k\rightarrow\infty}{\bf S}_{n_k}$ exists, and $|{\bf P}_{n_k}|\rightarrow \infty$, then  
$\lim_{k\rightarrow\infty}{\bf S}_{n_k} = \lim_{k\rightarrow\infty}{\bf P}_{n_k}$, so we will assume that any sequence $({\bf s}_{k})$ is periodic.

In this case  
Equation \ref{yuk} reduces to
\begin{equation}
\label{lesseryuk}
{\bf S}_{np}={\bf P}_{np}=\left(\prod_{m=0}^{p-1}\btau^{m}{\bf s}_{p-m}\right)
\prod_{l=0}^{n-1}\left(\prod_{i=0}^{p-1}
\btau^{lp+i}({\bf s}_{p-i})\right)
\end{equation}

Here ${\bf S}_{(n+1)p}$ is a prefix of ${\bf S}_{np}$ for $n\geq 1$, 
and 
if  $|{\bf S}_{np}|\rightarrow \infty$,
then
$\by:=\lim_{n\rightarrow \infty}{\bf S}_{np}$ exists,  and  by Proposition \ref{sigmainjectivity}  has at  least 2 $\sigma$-preimages. Note that $\by$ satisfies the equation 
\begin{equation}\label{theone}\left(\prod_{m=0}^{p-1}\btau^{m}{\bf s}_{p-m}\right)\,\btau^{p}({\bf y})= {\bf y}\,.\end{equation}
A similar argument applied to the sequence $({\bf S}_{(n+1)p+i})_{n}$ yields the other branch points.

 Conversely, if $\by$ is a non-$\btau$-fixed $k$-branch point, then
 there are $k$ letters $\{a_{1}^{n},a_{2}^{n},\ldots a_{k}^{n}\}$
 whose $\btau^{n}$ words share a maximal common proper suffix ${\bf
 T}_{n}$ and such that $\by\in [{\bf T}_{n}]$. Moving to a subsequence
 if necessary, we can assume that these $k$ letters are independent of
 $n$ and so if $A_{1}=\{a_{1}, \ldots a_{k}\}$, the sequence $({\bf
 T}_{n})$ is identical to the sequence $({\bf S}_{n})$ generated by
 $A_{1}$. The result follows.
\item

Statements 3  and 4  now follow from Part 2.
\endproof

\end{enumerate}
{\bf Examples}

\begin{enumerate} \setcounter{enumi}{3}

\item(Example 1) If $\btau(a)= aab$, and $\btau(b)=abb$, then the
  unique branch point $\by$ satisfies $b\btau (\by)= \by$.

\item(Example 2) If $\btau(a)=acb$, $\btau(b)=aba$ and $\btau(c)=aaa$,
  then only $\bu$ is a 2-branch point.

\item(Example 3) If $\btau(a)=acb$, $\btau(b)=aba$ and $\btau(c)=aca$,
  then $\{a,b\}$ is a suffix cycle, so $\bu$ is a 2-branch
  point. Letting ${A}_{1} =\{b,c\}$, we obtain the (only other)
  2-branch point $\by$ satisfying $a\btau(\by)=\by$, and $\by\neq
  \bu$.

\item If $\btau(a)=abc$, $\btau(b)=aacc$, and $\btau(c)=abcc$, then
  $\bu$ is $\sigma$-invertible, and if $A_{1}=\{b,c\}$, ${\bf
  s}_{2k+1}=cc$ and ${\bf s}_{2k}=c$ for $k\geq 1$, and no other
  choice of $A_{1}$ generates any other eventually different sequence
  $({\bf s}_{k})$.  Thus $(X_{\btau},\sigma)$ has two branch points,
  $\by$ and $\by^{*}$, satisfying $\by=c\btau(cc)\btau(\by)$ and
  $\by^{*}= cc\btau(c)\btau(\by^{*})$, and these two points are
  distinct.

\item 
If $\btau(a)=bbad,$ $\btau(b)=ab$, $\btau(c)=ad$, and $\btau(d) =
dac$, then $\btau$ has 3 fixed points, with $\btau^{\infty}(a)$ a
3-branch point, $\btau^{\infty}(b)$ a 2-branch point, and where
$\btau^{\infty}(d)$ is $\sigma$-invertible.  If $A_{1}= \{a,b\}$ then
the words $\{\btau^{n}(a), \btau^{n}(b)\}$ have $\btau^{n}(b)$ as a
maximal common suffix, and these words also converge to $\bu$.

\end{enumerate}

   We remark that the notion of a branch point has been studied, for
example in \cite{BDH}.  There the authors work with two dimensional
sequences and define the notion of a {\em right asymptotic} point. The
relationship is that if $\bx$ is a one-sided $k$-branch point, then
this leads to $k$ right asymptotic two sided points. The authors find
upper bounds on the number of right asymptotic points - being $2n-1$
if $|{\cal A}|=n$ and $\btau$ is left proper, and $2(n-l)+rl$ in
general, if $\btau$ has $l$ proper prefixes and $r$ proper suffixes
(Proposition 9 and Theorem 1 of \cite{BDH}), answering a question in \cite{lzpre}. We digress here to  show that if $2\leq i
\leq 2n-1$, there exists a left proper $\btau$ which has $i$ right
asymptotic points. A similar argument will work to create non-proper substitutions having a prescribed (appropriately constrained) number of branch points.

\begin{proposition}
\label{barge1}
Let $\btau$ be left proper on $n$ letters.
For any $m$-tuples of natural numbers 
$(k_{1},k_{2},\ldots , k_{m})$, $(j_{1},j_{2}, \ldots j_{m})$ where 
\begin{enumerate}
\item
$2 \leq k_{1} < k_{2} < \ldots < k_{m}\leq n$
\item $k_{1}+j_{1}(k_{1}-1)\leq n+1$, and 
for each $i\geq 2$, $k_{i}+ j_{i}(k_{i}-1)\leq n+1 - \sum_{l=1}^{i-1}k_{l}+j_{l}(k_{l}-1)$ and, 
\item
$\sum_{l=1}^{m} k_{l}+j_{l}(k_{l}-1)\leq n+m$,
\end{enumerate}
 there exists a left proper $\btau$ on $n$ letters, with $j_{i}$
 $k_{i}$-branch points, for $1\leq i \leq m$, and no other branch
 points.
\end{proposition}
{\bf Proof:} We first prove the special case $m=1$, $k_{1}=2$, and
$j_{1}\leq n-1$, which would generate a substitution with $j_{1}$ $2$-branch
points, and so $2j_{1}$ left asymptotic points. Letting ${\cal A}=
\{x_{1}, \ldots x_{n}\}$, choose some letter $a$ so that all $\btau$-substitution words start with $a$. Let $\btau(x_{1})$ have $x_{1}$ as suffix and 
let all other $\btau(x_{i})$'s have $x_{2}$ as suffix. As long as we
ensure later that both $x_{1}a,\, \, x_{2}a\,\, \in {\cal L}_{\btau}$,
then the $\tau$-fixed point is a $2$-branch point.  We define
$\btau(x_{i})$ by successively adding letters from the right hand
side, moving towards the left. Define the penultimate letter of
$\btau(x_{2})$ to be $x_{2}$, and the penultimate letter of all other
$\btau(x_{i})$ to be $x_{3}$. Then we have a second branch point
(satisfying the equation $x_{2}\btau({\bf y}) = \btau({\bf y})$). At
the $l$th stage (for $l\leq n-1$), fill in the $l$th column from the
right of $\btau(x_{i})$ by putting down an $x_{l}$ for $\btau(x_{l})$,
and $x_{l+1}$ elsewhere. In this way you get $n-1$ 2-branch points,
$n-2$ of them satisfying $x_{l}(\btau({\bf y})) = \btau({\bf y})$ for
$2\leq l \leq n-1$ (and the last $\btau$-fixed).  If $j_{1}<n-1$, you
stop the above procedure after having defined the $j_{1}$-st column
from the right. To ensure you get no more branch points (and also that
$\btau$ is injective, if required) start filling in columns with $x_{1}$ and $x_{2}$ 
(or $x_{i}$ and $x_{i+1}$, where $i\leq j_{1}$). Then any new branch points satisfy an already
defined equation - $x_{i}\btau(\by) = \by$ for some $i\leq j_{1}$.

In fact the above special case where $j_{1} = n-1$ is how one
maximises the number of branch, and left asymptotic, points. These
$2$-branch points can be packed in an optimal way. To obtain $j$
$k$-branch points, one needs $k+ k(j-1)$ letters, replicating the
above example, but at each stage when defining the $\btau$
substitution words, having $k$ different letters in each column that
you define. To ensure the branch points that you create are distinct,
the letters should change from column to column. When the desired
number of branch points have been obtained, start filling in new
columns with combinations of letters that have already been used to
create a previous branch point (for example use a subset of $x_{1},
x_{2}, \ldots x_{k}$ if they were used to create a $k$-branch point).

\endproof

Unilateral recognizability is less amenable than  bilateral
recognizability, which all aperiodic primitive substitutions possess,
and for whom the sequence of partitions $({\cal P}_{n})$ defined by
\[{\cal P}_{n}
:=
\{[\sigma^{k}(\btau^{n}(a))]: a \, \in {\cal A},\,\,\, 0 \leq k < |\tau^{n}(a)|\}
\]
is a nested sequence which spans the topology of $(X^{\mathbb
Z}_{\btau},\sigma)$ (see Proposition 14, \cite{dhs}). In the one sided
subshift $(X_{\btau}, \sigma)$, the sets in ${\cal P}_{n}$ need
not be disjoint. For example, if $\btau$ has a branch point $\by\neq
\bu$, then, for each $n$, there are letters $a$, $b$, and numbers
$n_{1}$ and $n_{2}$ such that $|\btau^{n}(a)|-n_{1} =
|\btau^{n}(b)|-n_{2}$ and
$\sigma^{n_{1}}(\btau^{n}(a))=\sigma^{n_{2}}(\btau^{n}(b))$, while
$\sigma^{n_{1}-1}(\btau^{n}(a))\neq\sigma^{n_{2}-1}(\btau^{n}(b))$. If
$\btau$ is such that all  branch points are fixed in $X_{\btau}$,
this does not happen.

\begin{proposition}
\label{clopenpartition}
Suppose that the  $k$ branch points $\{\bu_{1},\ldots \bu_{k}\}$ of $\btau$
are $\btau$-fixed.

Then
\begin{enumerate}
\item
For every $n>0$, ${\cal P}_{n}$ is a clopen partition, and the sequence of bases is decreasing. 
\item For each $n$, ${\cal P}_{n+1}$ is a refinement of ${\cal P}_{n}$: i.e.
every element of ${\cal P}_{n+1}$ is contained in an element of ${\cal
P}_{n}$.
\item
The intersection of the bases of ${\cal P}_{n}$ is $\{\bu_{1},\ldots \bu_{k}\}$, and the sequence $({\cal P }_{n} )$ spans the topology of $(X_{\btau},\sigma)$.

\end{enumerate}
\end{proposition}

{\bf Proof:}
\begin{enumerate}
\item
We show that ${\cal P}_{1}$ is clopen, First suppose that $({\bf
y}_{n})$ is a sequence in $\btau [a]$ and ${\bf y}_{n}\rightarrow
{\bf y}$. Since $\btau(X_{\btau})$ is closed, by Proposition
\ref{hosttopologyprop}, then ${\bf y}\, \in \btau(X_{\btau})$ If ${\bf
y}_{n}=\btau({\bf x}_{n})$ and ${\bf x}$ is a limit point of $({\bf
x}_{n})$, then ${\bf x}\,\in [a]$ and ${\bf y}=\btau({\bf x})\subset
[\btau(a)].$ Using Corollary \ref{almostuniquerep}, the sets
$\{\btau[a]: a\, \in {\cal A}\}$ are disjoint, and their union,
$\btau(X_{\btau})$, is a clopen set, so they are also.  If $0<k<|\btau (a)|$
and $({\bf y}_{n})$ is a sequence of elements in
$[\sigma^{k}(\btau(a))]$, then ${\bf y}_{n}=\sigma^{k}({\bf x}_{n})$
where ${\bf x}_{n}\,\in \btau([a])$. If ${\bf x}$ is a limit point of
$({\bf x}_{n})$ then it must also be in $\btau[a]$, and ${\bf y} =
\sigma^{k}({\bf x})\in \sigma^{k}(\btau[a])$. So the sets
$\sigma^{k}(\btau[a])$ are closed, and if the branch points are
all $\btau$-fixed, then these sets are also disjoint, and so clopen.

\item
 Take an element ${\bf p} \in {\cal P}_{n+1}$: we need to
show that $[\bp]\subset [\bp']$ where $[\bp'] \in {\cal P}_{n}$.
Suppose that 
$\bp =
\sigma^{k}(\btau^{n+1}(a))$, 
where $0<k< |\btau^{n+1}(a)|$. 
As in the proof of Proposition 14, in \cite{dhs}, write 
$\btau^{n+1}(a)= \btau^{n}(a_{1})\,\btau^{n}(a_{2})\, \ldots
\btau^{n}(a_{m})$ where $\btau(a)= a_{1}\, a_{2} \ldots a_{m}$. For
some $j$, $|\btau^{n}(a_{1}\ldots a_{j})|\leq k <
|\btau^{n}(a_{1}\ldots a_{j+1})|$, so let $l:=k -
|\btau^{n}(a_{1}\ldots a_{j})|$. Then
$\sigma^{k}(\btau^{n+1}([a]))\subset \sigma^{l}\btau^{n}[a']$, where 
$a'=a_{j+1}$. Let $\bp':= \sigma^{l}\btau^{n}[a']$.

 \item
 To see that $({\cal
P}_{n})$ span the topology of $(X_{\btau},\sigma)$, we imitate the proof of
Proposition 14, \cite{dhs}. Given $m$ positive, we show that for all  $n$
large,  each element ${\bf p}$ of
${\cal P}_{n}$ is contained in some ${\bf c}$ where ${\bf c}$ is a word of 
length $m+1$. Suppose that all $\btau$-substitution words start with
$l$. Write $l_{n} = |\btau^{n-1}(l)|$. Choose $n$ so that $l_{n}>m$
(minimality implies that $\lim_{n\rightarrow \infty}l_{n} = \infty$).

Fix $\bp \,\in {\cal P}_{n}$, and suppose that $\bp = \sigma^{k}(\btau^{n}(a))$ for some  $a\,\in {\cal A}$ and $k\,\in [0, |\btau^{n}(a)|)$. If $\by \,\in
[\bp]$ then there is some $\bx\,\in [a]$ with
$\by=\sigma^{k}\btau^{n}(\bx)$. Since $\btau$ is left proper,
${\btau}^{n}(a)\btau^{n-1}(l)$ is a prefix of ${\btau}^{n}(\bx)$, and
so $\by \in [\sigma^{k}{\btau}^{n}(a)\btau^{n-1}(l)], $ and this last
cylinder set has length greater than $m$. So, $\by$ starts with a
block of length $m$ depending only on $\bp$ and  not on $\by$.  \endproof
\end{enumerate}

\section{Bratteli Diagrams}\label{Bratteli}
A {\em stationary Bratteli diagram} ${\cal B}=({\cal V},{\cal E})$ is an infinite
directed graph with {\em vertex set} ${\cal V} = \bigsqcup_{n=
0}^{\infty} {V}_{n}$ and {\em edge set} ${\cal E} =
\bigsqcup_{n=1}^{\infty}{E}_{n}$, where  ${ V}_{n}=V$ for each $n\geq 1$, 
${E}_{n} = E$ for each $n\geq 1$, ${ V}_{0}= \{v_{0}\}$, and for each
$n\geq 1$ there is a {\em source map} $s:E_{n}\rightarrow V_{n-1}$ and
a {\em range map} $r:E_{n}\rightarrow V_{n+1}$.  If $|V_{n}|=d$, the
{\em incidence matrix} $F= (f_{i,j})_{1\leq i,\,j \leq d} $ is defined
by $f_{i,j}:= |\{ x \,\in E_{n}: s(x)=j, \,\,r(x)=i\}|$ for any $n\geq
2$. If $n=1$, the incidence matrix $F_{1}$ describing $E_{1}$ is a $d\times 1$ vector defined in a similar way.
We assume that
$s^{-1}(v)\neq \emptyset$ for each $v\in {\cal V}$ and $r^{-1}(v) \neq
\emptyset $ for all $v \, \in {\cal V}
\backslash v_{0}$. We will use $x, y\ldots$ when referring to edges,
and $a,b, \ldots$ when referring to vertices.
We say that ${\cal B}$ is {\em simple} if some power of $F$ has all positive entries.

A finite set of edges $\{x_{n+k}\}_{k=1}^{K}$, with
$s(x_{n+k+1})=r(x_{n+k})$ for $1\leq k\leq K-1$, is called a {\em path
} from $s(x_{n+1})$ to $r(x_{n+K})$; we can extend the maps $s$ and $r$ to paths 
by defining $s(\{x_{n+k}\}_{k=1}^{K}) = s(x_{n+1})$ and $r(\{x_{n+k}\}_{k=1}^{K}) = r(x_{n+K})$. Similarly an {\em infinite path}
in ${\cal B}$ is a sequence ${\bx}=(x_{n})^{\infty}_{n=1}$, with
$x_{n}\,\in {E}_{n}$ for $n\geq 1$, and $s(x_{n+1})=r(x_{n})$ for
$n\geq 1$.  
 The set of all infinite paths in ${\cal B}$
will be denoted $X_{\cal B}$ (a subset of $ \Pi_{n\geq 1} {E}_{n}
$), and  $X_{\cal B}$ is endowed with the topology induced from the product
topology on $\Pi_{n\geq 1} {E}_{n}$, so that $X_{\cal B}$ is a
compact metric space.
 If ${\bx}=
(x_{n})^{\infty}_{n=1}$ and ${\bx'}=(x_{n}')^{\infty}_{n=1}$ are
two elements in $X_{\cal B}$, we write ${\bx} \sim {\bx'}$ if there is some $N$ such that $x_{n}=x_{n}'$ for all $n\geq N$. It follows that $\sim$ is
an equivalence relation (the {\em tail equivalence relation}), and if ${\cal B}$ is simple, then each
equivalence class for $\sim$ is dense in $X_{\cal B}$, and $X_{\cal
B}$ has no isolated points, making it a Cantor space.

We will sometimes {\em telescope} ${\cal B}$ to ${\cal B'}$ - for our
purposes we need only define a restricted version of this
procedure. Let $k$ be given; then ${\cal B}^{(k)} =({\cal V}^{(k)},{\cal E}^{(k)})$
is a {\em telescoping}  of ${\cal B}= ({\cal V},{\cal E})$ if
${V}_{n}^{(k)}=V_{nk}$ (with the vertex $v \, \in {V}_{nk}$ labelled as $v
\, \in {V}_{n}^{(k)}$), and $F^{(k)}:=F^{k}$, so that the number of edges from
$v_{n}\, \in {V}_{n}^{(k)}$ to $v_{n+1}\, \in {\cal V}_{n+1}^{(k)}$ is the
number of paths from $v_{n}\, \in {V}_{nk}$ to $v_{n+1}\, \in
{V}_{(n+1)k}$. Conversely, we can perform a {\em splitting} by
introducing a new level between two consecutive levels ${V}_{n-1}$ and
${V}_{n}$, and new edge sets, so that if the resulting stationary
Bratteli diagram is ${\cal B'}=({\cal V}',{\cal E}')$, and we
telescope to every other level, we obtain the original ${\cal
B}=({\cal V},{\cal E})$. We will split diagrams only when a vertex in
$V_{1}$ has more than one edge connecting it to $v_{0}$. 
We consider
two Bratteli diagrams ${\cal B}$ and ${\cal B}'$ {\em equivalent} if
${\cal B'}$ can be obtained from ${\cal B}$ by relabeling,
telescoping and splitting, so that when we talk about a Bratteli diagram
we are talking about an equivalence class of diagrams.

\subsubsection{Ordering  $X_{\cal B}$}
 For each $n\geq 1$ and each $a \, \in {V_{n}}$, we put a linear order
 on $\{x\, \in E_{n}: r(x) = a\}$. In this article we assume that
 the order we put on $\{x\, \in E_{n}: r(x) = a\}$ depends only on
 $a$ and not $n$ if $n>1$ . We also call this ordering {\em stationary}. Note that we do not put any conditions on the ordering of $E_{1}$.
Elements of ${E}_n
 (a)$ will then be labelled $0, 1,2, \ldots$ according to their order,
 and we say ${\cal B}$ is {\em ordered }. We will also write $x<y$ if
 $x$ and  $y$ have the same range and $x$ has a smaller labelling than
 $y$. Every stationary ordering is generated by a substitution $\btau$
 on ${\cal A}$, an alphabet of size $|V|$: if the $i$th incoming edge
 to $a$ has source $b$, then define the $i$th letter of $\btau(a)$ to
 be $b$. Conversely to each substitution $\btau$ on ${\cal A}$, we can
 define the (stationary, ordered) {\em Bratteli diagram associated
 with $\btau$}, denoted ${\cal B }_{\btau}$: there is exactly one edge
 from each vertex in ${V}_{1}$ to $v_{0}$, let the transition matrix
 for $\btau$ be the incidence matrix $F$ for ${\cal B}$, and order
 ${\cal B}$ using the order of appearance of letters in the
 substitution words: if the $i$-th letter of $\btau(a)$ is $b$, then
 give one edge $x$ with $s(x)=b$ and $r(x)=a$ the ordering $i$. If
 ${\cal B}$ is the diagram defined by the transition matrix $M$ for
 $\btau$ then ${\cal B}^{(N)}$ is the diagram defined by the
 transition matrix $M^{N}$ for $\btau^{N}$.  We will use the
 substitution $\btau$ to denote the ordering on ${\cal B}$ and write
 ${\cal B}_{\btau}$  for the ordered diagram - in this case
 ${\cal B}^{(N)}_{\btau^{N}}$ is the ordered diagram associated to
 $\btau^{N}$.

For $m<n$, the order $\btau$ induces a partial ordering on paths from
${V}_{m}$ to ${V}_{n}$: the two paths ${\bx}=x_{m+1}\ldots x_{n}$ and
${\bx}'=x_{m+1}'\ldots x_{n}'$ with source in ${E}_{m}$ and common
range in ${E}_{n}$ are comparable, with ${\bx}<{\bx}'$ if there is
some $k\, \in [m+1,n]$ with $x_{k}<x_{k}'$ and $x_{j}=x_{j}'$ for
$k+1\leq j\leq n$.  If $1 \leq m<n$, then the number of paths in
${\cal B}$ from $b\,\in {\cal V}_{m} $ to $a\,\in {\cal V}_{n}$ is the
number of occurrences of $b$ in ${\btau}^{ (n-m)}(a).$ If $\btau$ is
primitive then there is a positive $k$ such that for any two letters
$a$ and $b$, there is at least one path from $ b\,\in {\cal V}_{n}$ to
$a\,\in {\cal V}_{n+k}$.

The ordering $\btau$ on ${\cal B}$ extends to an order on tail
equivalence classes in $X_{\cal B}$: if $\bx\sim \bx'$, we write
${\bx} < {\bx'}$ if an initial segment $x_{1}\ldots x_{N}<x_{1}'\ldots
x_{N}'$.  An infinite path is {\em maximal} ({\em minimal}) if all the
edges making up the path are maximal (minimal).  For any $\btau$ there
is a power $\btau^{N}$ such that all substitution words of $\btau^{N}$
have prefixes which generate right infinite $\btau$-fixed points, and
suffixes which generate left infinite $\btau$-fixed points.  Each one
of these prefixes generates a distinct minimal path in $X_{\cal B}$,
and each such suffix generates a distinct maximal path in $X_{\cal
B}$.  If $\{v_{1},v_{2},\ldots v_{n}\}$ ($ \{w_{1},w_{2},\ldots
w_{k}\}$) are these suffixes (prefixes), we call these the vertices
that determine the maximal (minimal) paths in $X_{\cal B}$. Note that the orderings $\btau$ and $\btau^{N}$ generate the same number of maximal and minimal paths in   
$X_{{\cal B}}$ and $X_{{\cal B}^{(N)}}$ respectively. If $\btau$
generates a minimal subshift, then $k\leq n$.  Note
that all maximal and minimal paths in $X_{{\cal B}}$ are
(shift) periodic, and those in $X_{{\cal B}^{(N)}}$ are
(shift) fixed.  Let $X_{\max}:=\{ {\bf M_{1}},
\, {\bf M_{2}},\ldots {\bf M_{n}}\}$ and  $X_{\min}=\{{\bf m_{1}},
{\bf m_{2}}, \ldots {\bf m_{k}}\}$  denote the sets of  maximal and minimal paths respectively  in 
 $X_{{\cal B}}$. They are both non-empty.

If ${\bx}= (x_{n})^{\infty}_{n=1}\, $ is not maximal, let $k$ be the
smallest integer such that $x_{k}$ is not a maximal edge, and let
$y_{k}$ be the successor of $x_{k}$ in $\{x: r(x)=r(x_{k})\}$. Then
the {\em successor} ${V}_{\btau} ({\bx}) $ is defined to be
${V}_{\btau}({\bx})= 0\ldots  0\, y_{k} x_{k+1} x_{k+2}\ldots $,
where $0\ldots 0$ is the unique minimal path 
ending at the source of $y_{k}$. Similarly, every non-minimal path has
a unique {\em predecessor}.  Simple ordered Bratteli diagrams which
have a unique minimal and maximal element are called {\em proper}, and
those with a unique minimal element ${\bf m}$ are called {\em
semi-proper}.  If ${\cal B}_{\btau}$ is
proper, then $V_{\btau}$ is a homeomorphism, and   if ${\cal B}_{\btau} $ is semi-proper, then $V_{\btau}$
can be extended to a continuous surjection on $X_{\cal B}$ by setting
$V_{\btau} ({\bf M}_{i} )= {\bf m}$ for all maximal elements ${\bf
M}_{i}$. If ${\cal B}_{\btau}$ is not semi-proper, the continuity of
$V_{\btau}$ turns out to depend entirely on ${\cal L}_{\btau}$: we
address this in the next  proposition. 
We
remark that a special case of this result is implicit in the formula on Line
12, Page 5 of \cite{lzpre}, when the number of maximal elements equals
the number of minimal elements in ${\cal B}_{\btau}$.

\begin{proposition} 
\label{closingproposition}
 Suppose that ${\cal B}= {\cal B}_{\btau}$ is a stationary ordered  Bratteli diagram, which has maximal elements $\{{\bf M}_{i}\}_{i=1}^{n}$ and minimal elements
 $\{{\bf m}_{i}\}_{i=1}^{k}$ where $k\leq n$, determined by the vertices  $ \{v_{1},v_{2},\ldots
 v_{n}\}$  and $ \{w_{1},w_{2},\ldots
 w_{k}\}$ respectively. 
 Then $V_{\btau}$ can be
 extended to a continuous mapping $V_{\btau}: X_{\cal B}\rightarrow
 X_{\cal B}$ if and only if for each $v_{i}\, \in \{v_{1},v_{2},\ldots
 v_{n}\}$ there exists a unique $w_{j}\,\in \{w_{1},w_{2},\ldots
 w_{k}\}$ such that $v_{i}w_{j}\,\in {\cal L}_{\btau}$.

\end{proposition}

{\bf Proof:}   
 We first prove the result for substitutions which
 generate maximal and minimal paths in $X_{\cal B}$ that are
 $\sigma$-fixed.
 In other words, suppose  that all $\tau$-substitution words have a prefix (suffix)
 which generates a right (left) infinite fixed point. Suppose also  that all words of length two in ${\cal L}_{\btau}$ appear in the $\btau$-substitution words. 

 First suppose that $V_{\btau}$ can be extended to a continuous map
 $V_{\btau}:X_{\cal B}\rightarrow X_{\cal B}$. 
If $f({\bf M}_{i})={\bf
 m}_{j}$, then for large $N$, elements in the cylinder set $[M_{i}^{N}
 e]$, with $e$ not a maximal edge,  are
 mapped, using continuity, into the cylinder set $[m_{j}^{N-1}]$.  If  $s(e)=s(e+1)=v$ 
 then $v_{i}w_{j}$
 is a subword of $ \btau^{2}(v)$,  and so in
 ${\cal L}_{\btau}$.
Suppose that $v_{i}w_{j'}\,\in{\cal L}_{\btau}$; we show that $j'=j$. 
Suppose  that $v_{i}w_{j'}$
 appears in a $\btau$-substitution word $\btau(v)$, then there exist
 edges $\alpha$, $\alpha +1$ with sources $v_{i}$ and $w_{j'}$ respectively
 and range $v$. For large $N$, take the cylinder set $[M_{i}^{N}\alpha]$; its image
 under $V_{\cal B}$ lies in $[m_{j'}^{N}\alpha+1]$ and  continuity
 means that $j=j'$.

Conversely suppose that for each $i$, there is a unique $j$ with
$v_{i}w_{j}\,\in {\cal L}_{\btau}$. Define $V_{\btau}({\bf M}_{i}) :=
{\bf m}_{j}$.  If ${V}_{\btau}$ is not continuous, then there are
distinct edges $\alpha$, $\beta$ with source $v_{i}$, but for large
$N$, $V_{\btau}(M_{i}^{N}\alpha)\, \in [w_{j'}^{N-1}] $ and $V_{\btau} (M_{i}^{N}\beta) \in [w_{j''}^{N-1}]$, where $j'\neq j''$, and one
of $j'$, $j''$, say $j'$, is distinct from $j$. If $r(\alpha)=v$, then
$v_{i}w_{j'} \, \in \btau^{2}(v)$, so that $v_{i}w_{j'}\,\in {\cal
L}_{\btau}$, a contradiction.

For a general substitution $\btau$, there is some power $\btau^{N}$
whose Bratteli diagram ${\cal B}^{(N)}$ satisfies all our initial
assumptions. Also, $(X_{\cal B}, V_{\btau})$ is continuous if and only if 
$(X_{{\cal B}^{(N)}}, V_{\btau^{N}})$ is continuous. Finally   the substitutions $\btau$ and
${\btau}^{N}$ generate the same sets $ \{v_{1},v_{2},\ldots v_{n}\}$
and $\{w_{1},w_{2},\ldots w_{k}\}$, and also ${\cal L}_{\btau}= {\cal L}_{\btau^{N}}$. The result follows.


\endproof

Let us call a stationary ordered diagram ${\cal B}_{\btau}$ {\em
closing} if it satisfies the conditions of Proposition
\ref{closingproposition}.  If ${\cal B}_{\btau}$
is closing, then $(X_{\cal B},V_{\btau})$ is called a {\em
Bratteli-Vershik} or {\em adic} system.  Since $V_{\btau}$ orbits are
equivalence classes for $\sim$, if $\btau$ is primitive, so that
${\cal B}_{\btau}$ is simple, then every $V_{\btau}$ orbit is dense.
Note that if there are multiple edges between $v_{0}$ and vertices in $V_{1}$, this has no bearing on the continuity of $V_{\btau}$: ie for any transition matrix $F_{1}$ and any 
order put on $E_{1}$, the resulting map $V_{\btau}$ (we abuse notation
here, as different $F_{1}$'s generate different $V_{\btau}$'s) is continuous if and only if the $V_{\btau}$ is continuous for the choice $F_{1}:= (1,1, \ldots 1)$.



The following is the specialised version of Theorem 4.6 in \cite{hps} that we will need:
\begin{proposition}
\label{rightadicrep}
Let $(X_{\btau},\sigma) $ be a substitution subshift where all $k$ 
branch points $\{\bu_{1},\bu_{2},\ldots \bu_{k}\}$
are $\btau$-fixed. If $V_{\btau}$ is continuous, 
 then $(X_{\btau},\sigma) 
\stackrel{\Phi}{\cong} (X_{\cal B},V_{\btau})$, with the points $\{\bu_{1},\bu_{2},\ldots \bu_{k}\}$ mapped to the $k$ minimal elements of $X_{\cal B}$.
\end{proposition}

{\bf Proof:} This follows from Proposition \ref{clopenpartition}, and standard constructions of $\Phi$, as for example in the proof of Theorem 4.6 in \cite{hps}.

\section{Quasi-invertible substitutions}\label{quasi_invertible}

In this section we obtain adic representations for any
quasi-invertible substitution $\btau$. Here we are in the situation where we
do not have to worry about the continuity of $V_{\btau}$, so that to
apply Proposition \ref{rightadicrep} we need to find some $\btau'$ with $(X_{\btau},\sigma)\cong  (X_{\btau'}^{\mathbb N},\sigma)$, and where the branch point of $\btau'$ is $\btau'$-fixed.

\begin{lemma}
\label{nonfixedlemma2a}
Suppose that  $\btau$ is defined as  $\btau (a_{i})= {\bf p}_{i} x_{i} {\bf s}_{1}$, where
$\{x_{1},\ldots, x_{n}\}$ has at least two distinct elements, and where if ${ A}_{1}= {\cal A}$, then  
${\bf s}_{k}= {\bf s}_{1}$ for $k\geq 1$.  Define
the substitution ${\btau}^{*}$ as ${\btau}^{*}(a_{i}) = {\bf
s}_{1}{\bf p}_{i}x_{i}$.
\begin{enumerate}
\item
For each $n\geq 0$, and each $\bw \, \in {\cal L}_{\btau}$,
${\btau}^{*}({\btau}^{n}(\bw)){\bf s}_{1} = {\bf s}_{1} \btau^{n+1}(\bw)\,.$
\item
 If ${\bf S}_{n}$ is as in Lemma \ref{nonfixedlemma1a}, then for each $n \geq 2$,
${\bf S}_{n}= \btau^{*}({\bf S}_{n-1})\,{\bf s}_{1}\, .$
\item
${\bf S}_{n}\subset ({\btau}^{*})^{n-1}({\bf s}_{1})$
for each $n> 1$.
\end{enumerate}

\end{lemma}

{\bf Proof}:
\begin{enumerate}
\item
Note that 
$\btau^{*}(a)\, {\bf s}_{1}= {\bf s}_{1}{\bf p}_{a}x_{a}{\bf s}_{1} = {\bf s}_{1}{\btau}(a)\, $ for any $a\in {\cal A}$.  The general proof for longer words follows by concatenation.
\item
 First by Lemma \ref{nonfixedlemma1a} and Part 1 of this lemma, 
${\bf S}_{2}= {\bf s}_{1}\btau ({\bf s}_{1})= \btau^{*}({\bf s}_{1}){\bf s}_{1}$. Assuming that ${\bf S}_{n}= \btau^{*}({\bf S}_{n-1}){\bf s}_{1}$, we have 
\[{\bf S}_{n+1} = {\bf S}_{n} \btau^{n}({\bf s}_{1})
\stackrel{IH}{=} \btau^{*}({\bf S}_{n-1}){\bf s}_{1} \btau^{n}({\bf s}_{1})
\stackrel{Part 1}{=} 
\btau^{*}({\bf S}_{n-1}){\btau}^{*}(\btau^{n-1}({\bf s}_{1})){\bf s}_{1}
= \btau^{*}({\bf S}_{n})  {\bf s}_{1}.\]
\item
The case $n=2$ is clear.
Assuming ${\bf S}_{n}\subset ({\btau}^{*})^{n-1}({\bf s}_{1}),$
\[{\bf S}_{n+1} \stackrel{Part 2}{=} {\btau}^{*}({\bf S}_{n}) {\bf s}_{1}
\stackrel{IH}{\subset} {\btau}^{*}(({\btau}^{*})^{n-1})({\bf s}_{1})) {\bf s}_{1}
= ({\btau}^{*})^{n}({\bf s}_{1}) \, {\bf s}_{1}\subset ({\btau}^{*})^{n}({\bf s}_{1})\,.\]
\endproof\end{enumerate}

\begin{corollary}
\label{nonfixedcorollary}

Let $\btau$ be  left proper and quasi-invertible,    
with 
branch point $\by\neq \bu$. Then some power of $\btau$ is right proper, and there 
exists a left proper, quasi-invertible substitution $\btau^{*}$ with  
$\btau^{*}(\by) = \by$, and $(X_{\btau},\sigma)=(X_{\btau^{*}}, \sigma)$.
\end{corollary}

{\bf Proof}: Since the fixed point $\bu$ is $\sigma$-invertible, Proposition  \ref{inverseimagesofu} tells us that some $\btau$ (or some power of $\btau$) is right proper.
If $\by$ is the  branch point for $\btau$,
then Lemma \ref{nonfixedlemma1a} tells us that there exists an ${\bf s}$ such that if $A$ is any subset of  ${\cal A}$, it  generates an eventually fixed sequence $({\bf s}_{k})$ where  ${\bf s}_{k}= {\bf s}$ for $k\geq n_{A}$.
 Since there are finitely many such sequences we can assume, taking a power if necessary,  that $n_{A}=1$.  Now Lemma \ref{nonfixedlemma2a} can be applied, taking 
$\btau^{*}$ as defined. The branch point is  $\by:=\lim_{n\rightarrow \infty}{\bf S}_{n}$; that it
exists and has at least two $\sigma$-preimages follows from Lemma
\ref{nonfixedlemma1a}. That it is ${\btau^{*}}$-fixed follows from
Part 3 of Lemma \ref{nonfixedlemma2a}. Thus $X_{\btau^{*}}\subset
X_{\btau}$, and by minimality, this inclusion is an equality.
\endproof

The connection between Bratteli-Vershik systems and one sided left
proper  substitutions is given by the next result. It is the
appropriate generalization of \cite[Prop 20]{dhs} to one sided
substitution systems.

\begin{theorem}
\label{quasiinvertibleleftproperareadic}

If $\btau$ is left proper and quasi-invertible, with fixed
point $\bu$, there exists a semi-proper, stationary Bratteli diagram
${\cal B}$ which is semi-proper, and stationary, such that $(X_{\cal
B}, V_{\cal B})$ is topologically conjugate to $(X_{\btau}, \sigma)$.

\end{theorem}

{\bf Proof}: If $\bu$ is the branch point,  the sequence of partitions ${\cal
P}_{n}$ defined in Section \ref{suffix} are a refining sequence of
partitions, by Proposition \ref{clopenpartition}, which generate the
topology of $(X_{\btau}, \sigma)$. Proposition \ref{rightadicrep} tells us that 
$(X_{\btau},\sigma) \cong (X_{\cal B},V_{\btau})$.
If $\btau$ is left proper and the branch point $\by\neq \bu$, then by
Lemma \ref{inverseimagesofu}, $\btau$ (or some power of $\btau$) has
to be right proper.  Work with the substitution $\btau^{*}$ in
Corollary \ref{nonfixedcorollary}. As above, $(X_{\btau^{*}},\sigma)$
is conjugate to $(X_{{\cal B}^{*}}, V_{{\btau}^{*}})$, and since
$(X_{\btau},\sigma)=(X_{\btau^{*}},\sigma)$, the result
follows.\endproof

\subsection{Return words and induced substitutions \label{inducedsubstitutions}}
We now introduce concepts needed to extend Theorem
\ref{quasiinvertibleleftproperareadic} to non-left proper
substitutions.
Suppose that $T:X\rightarrow X$ is a minimal continuous transformation
 with $X$ a Cantor space - henceforth called a Cantor system.  If
 $U\subset X$ is a clopen set, the {\em system $(U,T_{U})$ induced by
 $(X,T)$ on $U$}, is defined by $T_{U}(\bx) = T^{n}(\bx)$, where $n$
 is the least positive natural number such that $T^{n}(x) \, \in
 U$. $T_{U} $ is well defined since $T$ is minimal, and since $X$ is
 compact, $n$ can only take a finite number of values.  The Cantor
 systems $(X_{1},T_{1})$ and $(X_{2}, T_{2})$ are (topologically) {\em
 Kakutani equivalent} (\cite{Pe} is first reference of this) if there
 exist clopen sets $U_{i}\subset X_{i}$ such that the respective
 induced systems are isomorphic. Topological Kakutani equivalence is
 much more stringent than measurable Kakutani equivalence - for
 example, all rank one transformations are measurably Kakutani
 equivalent, but it can be shown, for example, that the {\em Chacon}
 substitution $\btau(0)=0010,$ $\btau(1)=1$ is Kakutani equivalent
 only to the substitutions $\btau_{1}(0)=001^{k}0$, $\btau(1)=1$, or
 $\btau_{1}(0)=01^{k}00$, $\btau_{2}(1)$.

   Conversely, given an induced  system $(U,\sigma_{U})$, where $(X,\sigma)$ 
is minimal, for some finite $H$, we can define the {\em height function} $h:U\rightarrow \{1,\ldots H\}$ given by $h(\bx)=k$ if and only if 
$\sigma^{i}(\bx )\not\in U$ for $1\leq i\leq k-1$ but $\sigma^{k}(\bx)\in U$.
 For each $k$ in $\{1,\ldots H\}$, and $0\leq k-1$, let $U_{k}^{i}=
\sigma^{i}(h^{-1}(\{k\}))$. Since $(X,\sigma)$ is minimal,
$X=\cup_{k=1}^{H}\cup_{j=0}^{k-1}U_{k}^{j}$.

\begin{lemma}
\label{clopenlemma}
If $(X_{\btau},\sigma)$ is quasi-invertible, with branch point in the
clopen set $U$, then the sets $\{\{U_{k}^{j}\}: 1\leq k\leq H, 0\leq
j\leq k-1\}$ form a clopen partition of $X_{\btau}$.
\end{lemma}

{\bf Proof:} It is straightforward that each $U_{k}^{0}$ is clopen, by
their definition. If $({\by}_{n}) \subset U_{k}^{i}$ and ${\by}_{n}
\rightarrow {\by}$, then for each $n$ there is some ${\bf x}_{n}\,\in
U_{k}^{0}$ with $\sigma^{k}({\bf x}_{n}) = {\bf y}_{n}$. If ${\bf x}$
is a limit point of the points ${\bf x}_{n}$, then ${\bf x}\,\in
U_{k}^{0}$ and $\by=\sigma^{i}(\bx)$. Thus $U_{k}^{i}$ is closed. If the branch point in $X_{\btau}$ is in $U$, these sets are also disjoint. The result follows.
\endproof

 If we can partition $X_{\btau}$ as in Lemma \ref{clopenlemma},  we say that $(X_{\btau},\sigma)$ is a {\em primitive} of $(U,\sigma_{U})$.
If $(U,T_{U})$ has an appropriate adic representation, it is now
possible to extend this representation to $(X_{\btau},T)$, as in the proof of
Theorem 3.8, in \cite{gps}:

\begin{theorem}
\label{gpstower}
Let $(X_{\cal B}, {V}_{\cal B})$ be an adic system, where $\cal B$ is
semi proper, simple, and with M maximal elements. Suppose that the
aperiodic M-quasi-invertible Cantor system $(Y,T)$, with branch point
$\by$ has an induced system $(Y',S)$ where $\by\, \in Y'$. If $(Y',S)$
is isomorphic to $(X_{\cal B}, V_{\cal B})$, then $(Y,T)$ is
isomorphic to $(X_{\cal B'}, V_{\cal B'})$ where ${\cal B'}$ is
obtained from ${\cal B}$ by adding or removing a finite number of
edges to ${\cal E}_{1}$, and changing the ordering on the affected
vertices and edges.  \endproof\end{theorem}

Let $\btau$ be a (not necessarily left proper)  substitution with fixed point $\bu$, and suppose that $u_{0}=a$. A {\em return word to $a$} in $\bu$ is a word
${\bw}$ such that
\begin{enumerate}
\item
$u_{0}$ is a prefix of ${\bw}$;
\item
There is no other occurrence of $u_{0}$ in ${\bw}$;
\item
${\bw}u_{0}\, \in {\cal L}_{\btau}$.
\end{enumerate}

This definition is a special case of the definition of a return word in \cite{dua}.
The set of return words ${\cal R}$ is finite, since $\bu$ is almost
periodic.  This notion is the one-sided generalization of two sided
return words in \cite{dhs}, and we use their notation here.
The set of return words ${\cal R}$ is finite, since $\bu$ is almost
periodic.  This notion is the one-sided generalization of two sided
return words in \cite{dhs}, and we use their notation here.  Ordering
${\cal R}$ according to the order of appearance of a return word in
$\bu$, we have a bijection $\psi$ from $R:=\{1,2,\ldots |{\cal R}|\}$
to ${\cal R}$. The $\btau$-fixed point $\bu$, and so every element in
$X_{\btau}$, is a concatenation of return words, with possibly a suffix of a
return word as a prefix. Extend $\psi$ by concatenation to
$\psi:R^{\mathbb N}\rightarrow [u_{0}] \subset {\cal A}^{\mathbb
N}$. Since any $\by\,\in \psi(R^{\mathbb N})$ can only be partitioned
in one way, using Property 2 of the definition of a return word,
$\psi$ is injective.

 Let ${\cal D}(\bu) \in R^{\mathbb N}$ be
the unique sequence such that $\psi({\cal D}(\bu)) = \bu$.
Let $(Y,\sigma)$ be the subshift spanned by ${\cal D}(\bu) \in R^{\mathbb N}$.
\begin{lemma}
\label{kakutanilemma}
$(Y,\sigma)$ is isomorphic to the system induced by $(X_{\btau},\sigma)$ on
$[u_{0}]$ via the map $\psi$.
\end{lemma}

{\bf Proof:} For each $\by\, \in Y$, $\psi(\sigma \by) =
\sigma^{k}(\psi (\by)$, where $k$ is the length of the word in ${\cal
R}$ corresponding to $y_{0}$. So $\psi(\sigma^{n}(\by))\, \in [u_{0}]$
for all natural $n$, thus $\psi (Y) \subset [u_{o}]$. If $\bx \, \in
[u_{0}]$, then $\sigma^{n_k}(\bu) \rightarrow \bx$ for some
$n_{k}\rightarrow \infty$. For large $k$, $\sigma^{n_k}(\bu)\in
[u_{0}]$, so $\sigma^{n_k}(\bu) = \psi (\sigma^{j_k} ({\cal D}(\bu))\,
\in \psi (Y)$ for some $j_{k}\rightarrow \infty $. Hence $\bx\, \in
\psi(Y)$. It follows that $\psi \circ \sigma = \sigma_{[u_0]}\circ
\psi$ for $y\, \in Y$.  \endproof

 The substitution $\btau$ can be used to define a substitution on $R$:
 if ${\bw}$ is a return word, then ${\bw}$ starts with $u_{0}$, and
 since $\btau(u_{0})$ starts with $u_{0}$, so does $\btau
 ({\bw})$. Also, ${\bw}u_{0}$ is a word in $\bu$, so $\btau ({\bf
 w}u_{0})= \btau({\bw})u_{0}\ldots$ is also a word in $\bu$. Hence
 $\btau({\bw})u_{0}...$ is a concatenation of unique return words, and
 a prefix of some return word, so that $\btau ({\bw})$ is a unique
 concatenation of return words. Thus if $j\, \in R$ corresponds to
 ${\bw}\,\in {\cal R}$, and $\btau ({\bw})= {\bw}_{1}{\bw}_{2}\ldots
 {\bw}_{r}$, define $\btau_{1}(j): = i_{1}i_{2}\ldots i_{r}$ where
 $i_{k} \, \in R$ corresponds to ${\bw}_{k}\, \in {\cal R}$.

\begin{lemma}
\label{towerovertau1}
Suppose $\btau$ have generating fixed point $\bu$. Then
\begin{enumerate}
\item 
$\btau_{1}$ (or some power of $\btau_{1}$) is left proper, primitive,aperiodic  injective, with 
${\cal D}(\bu)$ its  fixed point;
\item
$\btau_{1}$ is recognizable; and
\item
 If $\btau$ is $M$-quasi-invertible, then $\btau_{1}$ is
$M$-quasi-invertible. Conversely, if $\btau$ is quasi-invertible, and
$\btau_{1}$ is $M$-quasi-invertible, then $\btau$ is
$M$-quasi-invertible.
\end{enumerate}

\end{lemma}

{\bf Proof:}
\begin{enumerate}
\item
 Note that if $1\,\in R$ corresponds to ${\bw}$ in ${\cal R}$, then
${\bw}u_{0}$ is a prefix of $\bu$, and so is $\btau^{n}(u_{0})$ for
each $n$. Choose $n$ large enough so that
$|\btau^{n}(u_{0})|>|{\bw}|$, so that ${\bw}u_{0}$ is a prefix of
$\btau^{n}(u_{0})$. If $j\,\in R$, and ${\bw'}\, \in {\cal R}$
corresponds to $j$, then ${\bw'}$ begins with $u_{0}$, so that
$\btau^{n}(u_{0})$ and so ${\bw}$ is a prefix of
$\btau^{n}({\bw'})$. Hence 1 is the first letter of $\btau_{1}(j)$.

Since $\psi\circ\btau_{1}=\btau\circ\psi$, then
$\btau_{1}(D({\bu}))=D({\bu}))$. Since $\btau$ is primitive, so is $\btau_{1}$'
If $\btau_{1}(i)=\btau_{1}(j)$, then either $\btau
(\alpha)=\btau(\beta)$ for some $\alpha\neq \beta$, contradicting
injectivity of $\btau$, or the conditions in Theorem
\ref{mosserecognizability} are satisfied, contradicting
recognizability of $\btau$.
\item
 Suppose that $\btau_{1}$ is not recognizable. Then by Theorem
 \ref{mosserecognizability}, there exist letters $i$ and $j$ in $R$
 and some $\by \, \in Y$ with $\btau_{1}(i)$ a proper suffix of
 $\btau_{1}(j)$, and $\btau_{1}(i)\by$ and $\btau_{1}(j)\by$ appearing
 in $Y$ with the same 1-cutting of $\by$. Suppose that $\psi(i) =
 {\bw}= w_{1}\ldots w_{l}$ and $\psi(j)={\bw}'=w_{1}'\ldots
 w_{l}'$. Find the smallest $k$ in $\{0\ldots l-1\}$ such that
 $w_{l-k}=w'_{l'-k}$. Such a $k$ exists since otherwise $\bw$ is a
 suffix of $\bw'$, and since $\bw\neq \bw'$, then $w_{1}=u_{0}$ occurs
 in $w_{2}'\ldots w_{l}'$, a contradiction to Part 2 of the definition of a return word.
Now as $\btau_{1}(i)$ is a proper suffix of $\btau_{1}(j)$, there exist words $\bv$ and $\bv'$ in $R^{+}$, with $\bv$ a suffix of $\bv'$, $\bv $ a suffix of $\btau_{1}(i)$, $\bv'$ a suffix of $\btau_{1}(j)$, and words $\bp$, $\bp'$ and $\bs \, \in {\cal A}^{+}$, with 
$\psi(\bv)=\bp\btau(w_{l-k})\bs$ and $\psi(\bv')=\bp'\btau(w_{l'-k}')\bs$. If 
$|\btau(w_{l-k})|< |\btau(w_{l'-k}')|$, then $\btau(w_{l-k})$ is a proper suffix of $\btau(w_{l'-k}')$; and if $|\btau(w_{l-k})|> |\btau(w_{l'-k}')|$, then $\btau(w_{l'-k}')$ is a proper suffix of $\btau(w_{l-k})$. Finally if $|\btau(w_{l-k})|= |\btau(w_{l'-k}')|$, then $\btau(w_{l-k})=\btau(w_{l'-k}')$, which is not possible as $\btau$ is injective on letters.

Thus (without loss of generality) $\btau(w_{l-k})$ is a proper subword of $\btau(w_{l'-k}')$, and if $\by'=\btau(w_{l-k+1}\ldots w_{l})\psi(\by)$, then $\btau(w_{l-k})\by'$ and $\btau(w_{l'-k}')\by'$ appear in $X_{\btau}$ with the same 1-cutting, a contradiction to the recognizability of $\btau$.
\item
 Note that if $\bx\in X_{\btau}$ has distinct $\sigma$-preimages
$\bx_{1}, \,\bx_{2}, \ldots \bx_{M}$, then in $[u_{0}]$, there exist
$\by_{1},\, \by_{2}, \ldots \by_{M}$ with
$\by_{i}=\sigma^{-t_{i}}\bx_{i}$, where $-t_{i}$ is the first time a
preimage of $\bx_{i}$ is in $[u_{0}]$, - that this point is
well-defined follows from the quasi-invertibility and the aperiodicity
of $(X_{\btau},\sigma)$. Let $\by=\sigma_{[u_{0}]}(\bx)$ if
$\bx\not\in [u_{0}]$, and $\bx=\by$ otherwise. Then $\by$ has $M$
distinct $\sigma_{[u_{0}]}$-preimages, $\by_{1},\ldots, \by_{M}$.
Since $\bx$ is unique, so is $\by$.  Now use the isomorphism in Lemma
\ref{kakutanilemma} to transfer this information to $(Y,\sigma)$.

Conversely, if $(X_{\btau}, \sigma)$ is $N$-quasi-invertible and $\bx$
 has $N$ $\sigma$-preimages, $\{\bx_{1},\ldots \bx_{N}\}$, then if $t$
 is the first re-entry time of $\bx$ into $[u_{0}]$, $\sigma^{t}(\bx)$
 has $N$ $\sigma_{[u_0]}$-preimages $\{\sigma^{-t_1}(\bx_{1}),\ldots
 \sigma^{-t_N}\bx_{N}\}$.  Thus $\psi^{-1}(\sigma^{t}(\bx))$ has $N$
 $\sigma$-preimages, and as $(X_{\btau_1},\sigma)$ is
 $M$-quasi-invertible, $N=M$. \endproof
\end{enumerate}

\begin{theorem}
\label{quasiinvertiblearealmostadic}

Suppose that $\btau$  is  $M$-quasi-invertible, with a generating fixed point $\bu$, and branch point $\by$.
Then 
 there
exists an ordered semi-proper, stationary Bratteli diagram ${\cal B}_{\btau^{*}}$ such that $(X_{\cal B}, V_{\btau^{*}})$ is topologically
conjugate to $(X_{\btau}, \sigma)$.

\end{theorem}

{\bf Proof:} Using the fixed point $\bu$, and working with the induced
system $([u_{0}],\sigma_{[u_{0}]})$, we use Lemma \ref{towerovertau1}
to find $\btau_{1}$ such that $([u_{0}],\sigma_{[u_{0}]})$ is conjugate to
$(X_{\btau_{1}},\sigma)$.  By Part 3 of Lemma \ref{towerovertau1},
$\btau_{1}$ is $M$-quasi-invertible.  Now using Corollary
\ref{nonfixedcorollary}, there is some $\btau^{*}$ left proper,
and quasi-invertible, with its fixed point ${\bf u}^{*}$ as
branch point, such that $(X_{\btau_{1}},\sigma)= (X_{\btau^{*}},
\sigma)$ (If $\btau_{1}$'s branch point is fixed, let $\btau_{1}=
\btau^{*}$). By Theorem \ref{quasiinvertibleleftproperareadic} $(X_{{\btau}^{*}}, \sigma)$,
and so $(X_{\btau_{1}},\sigma)$, is isomorphic to $(X_{{\cal B}_{\btau^{*}}},
V_{\btau^{*}})$. If $\by\,\in [u_{0}],$ then by Lemma
\ref{clopenlemma}, $(X_{\btau},\sigma)$ is a primitive of
$([u_{0}],\sigma_{[u_{0}]})$, and $\by$ is mapped to the minimal
element of $(X_{{\cal B}_{\btau^{*}}}, V_{{\btau}^{*}})$.  Using Theorem
\ref{gpstower}, $(X_{\btau},\sigma)$ is isomorphic to $(X_{\cal B },
V_{\btau^{*}})$ where ${\cal B}$ is obtained from ${\cal B}^{*}$ by the
addition of a finite number of edges to $ E_{1}$, and $V_{\btau^{*}}$ is the corresponding adic map.

If $\by\not\in [u_{0}]$, write $(X_{\btau},\sigma)$ as a primitive over a
sufficiently small cylinder set $V$ containing $\by$, as in Lemma
\ref{clopenlemma}, and such that if a level of this partition ${\cal
P}$ intersects $[u_{0}]$, then it is contained in $[u_{0}]$: thus
${\cal P}$ has $[u_{0}]$ split up as a union of cylinder sets defined
by words of length $N$. We make $V$ small enough so that no element in
$V$ returns to $V$ before it passes through $[u_{0}]$.  Telescope  just the
first $M$ levels of the Bratteli diagram ${\cal B}_{\btau^{*}}$, so that the
elements in ${\cal P}$ which are subsets of $[u_{0}]$ are represented
by edges from $v_{0}$ to the first level. Now add levels in ${
E}_{1}$ according to how the cylinder sets in $[u_{0}]$ appear in the
clopen partition of $X_{\btau}$. For example, If $C$ and $C'$ are cylinder
sets contained in $[u_{0}]$ and appearing as $C=V_{k}^{i}$ and
$C=V_{k}^{i'}$ with $i<i'$ and no set $V_{k}^{i+1},\ldots
V_{k}^{i'-1}$ intersects $[u_{0}]$, then in the telescoped Bratteli
diagram, the edge $e'$ corresponding to $C'$ is the successor of the
edge $e$ corresponding to $C$; now insert $i'-i-1$ new edges between
$e$ and $e'$.  The fact that elements in $V$ have to pass through
$[u_{0}]$ before returning to $V$ means that you can add enough levels
to the telescoped stationary Bratteli diagram (which is still defined by $\btau^{*}$) so that the new Bratteli diagram
${\cal B}_{\btau^{*}}$ generates an adic system $(X_{{\cal B}^{*}},V_{\btau^{*}})$
isomorphic to the original $(X_{\btau},\sigma)$.

\endproof

Lemma 9 in \cite{f} is used in the two sided version of this previous
result, to show that (aperiodic, primitive) substitution subshifts are
conjugate to a stationary adic system where there are only single
edges from the vertex $v_{0}$ to any vertex in ${\cal
V}_{1}$. Although the statement of this lemma is also true for
semi-proper Bratelli diagrams, it is not clear that the resulting ordered
stationary diagram corresponds to a recognizable, or injective
substitution. See for example, Figure 4 in \cite{dhs}.

{\bf Examples:}
\begin{enumerate}\setcounter{enumi}{8}
\item
If $\btau$ is  suffix permutative, with one
fixed point $\bu$, then it is quasi-invertible, and
$(X_{\btau},\sigma)$ is a primitive over $(X_{\btau_{1}}, \sigma)$, as
$\phi(\bu) \, \in X_{\btau_{1}}$.

\item
 If $\btau(a)=aac$, $\btau(b)=bcc$, and $\btau(c)=abc$, then
$(X_{\btau},\sigma)$ is 3-quasi-invertible, with the non-invertible
element $\by$ satisfying $ c\btau(y)=y$. Since
all $\btau$-words have a proper common suffix, Corollary
\ref{nonfixedcorollary} tells us that
$(X_{\btau},\sigma)=(X_{\btau^{*}}, \sigma)$ where
$\btau^{*}(\by)=\by$. Now apply  Theorem \ref{quasiinvertibleleftproperareadic}.

\item \label{adic}

If $\btau(a)=aac$, $\btau(b)=bcc$, $\btau(c)=adbc$ and
$\btau(d)=adbd$, then since $da\not\in {\cal L}_{\btau}$ and
$cb\,\not\in {\cal L}_{\btau}$, both fixed points are
$\sigma$-invertible. Inspection of maximal proper common suffixes of
the $\btau$-substitution words of letters in any subset ${\cal A}^{*}$
of ${\cal A}$ leads to the existence of only one branch point $\by$
satisfying $c\btau(\by)=\by$. Lemma 10 cannot directly apply because
the family of $\btau$-substitution words does not have a common
suffix.

If we consider instead the induced system $([c], \sigma_{[c]})$,
$(X_{\btau}, \sigma)$ is a primitive over $([c],
\sigma_{[c]})$. Define $\btau^{*}(a)=caa$, $ \btau^{*}(b)=c$,
$\btau^{*}(c)=cadb$ and $\btau^{*}(d)=cadbdb$. One can now prove,
similarly to the proof of Lemma \ref{nonfixedlemma2a} to show that
$\btau^{*}(\by)=\by$. Thus $([c],\sigma_{[c]})$ is itself a (left
proper, quasi-invertible) substitution subshift, and so
has a stationary adic representation $(X_{\cal B^{*}},V_{\cal
B^{*}})$. Now Theorem \ref{gpstower} applies, so that
$(X_{\btau},\sigma)$ also has a stationary adic representation. In
this example it is relatively straightforward to find the right
$\btau^{*}$; making the adic representation of $(X_{\btau},\sigma)$
more straightforward than if we had followed the proof of Theorem
\ref{quasiinvertiblearealmostadic};
in general though it is not clear how a branch point
$\by$ can be seen as the fixed point of some $\btau^{*}$ which is obtained directly from the definition of $\btau$.

\item\label{chacon}
{\em Minimal rank one subshifts } are defined by substitutions on
${\cal A}=\{0,1\}$ of the form
\[\btau(0)= 0^{n_{1}}1^{m_{1}}0^{n_{2}}1^{m_{2}}\ldots
0^{n_{k-1}}1^{m_{k-1}}0^{n_{k}} \mbox{ and }\btau (1)=1,\] where
$k<\infty$, and $n_{1}$ and $n_{k}$ are positive. The latter condition ensures that the
resulting substitution subshift is minimal, with generating fixed
point $\bu:=\lim_{n\rightarrow \infty} \btau^{n}(0).$ These systems
are equivalent to rank-one systems defined by `cutting and stacking'
where there are a bounded number of cuts and spacers added, and the
same cutting-and-stacking rule is obeyed at each stage. A
comprehensive exposition of rank one systems is given in
\cite{fe}. 
Suppose $n_{i}\neq n_{j}$ for some $i,j$, so that $\btau$ is
aperiodic. Then $(X_{\btau},\sigma)$ is quasi-invertible if and only
if $m_{1}=m_{2}= \ldots = m_{k-1}=m>0$. For, the generating fixed
point $\bu$ is a branch point, and the fixed point $1^{m_i}\bu$ is another branch
point if and only 
$m_{i}<m_{j}$ for some $j$.

If $\btau(0)= 0^{n_{1}}1^{m}0^{n_{2}}1^{m}\ldots
0^{n_{k-1}}1^{m}0^{n_{k}} $, then there are two return words to 0,
$\bw_{1}=0$ and $\bw_{2}=01^{m}$. In this case $\btau_{1}$ is defined
on $\{1,2\}$ and $\btau_{i}=\bw i$ where $\bw =
1^{n_{1}-1}21^{n_{2}-1}2\ldots 1^{n_{k}-1}$. By Lemma
\ref{towerovertau1}, $\btau_{1}$ is left proper and
2-quasi-invertible. Thus $(X_{\btau},\sigma)$ is a primitive of
$([0],\sigma_{[0]})$ and by Theorem
\ref{quasiinvertiblearealmostadic}, $\btau$ has a semi proper adic
representation by adding $m-1$ edges from $v_{0}$ to the vertex in
${\cal V}_{1}$ corresponding to the letter `1' in the Bratteli
representation of $\btau_{1}$. For $1\leq k\leq M$, the $\btau$-fixed
points $1^{k}\bu$ live on (different) levels of the tower over
$(X_{\btau},\sigma)$.

If for some $i,j$, $m_{i}\neq m_{j}$, then $\btau_{1}$ is still of the
form $\btau_{1}(i)=\bw i$, where ${\bw}= 1^{n_{1}-1}2 1^{n_{2}-1}3
\ldots 1^{n_{k-1}-1}k1^{n_{k}-1}$, so that $\btau_{1}$ is
M-quasi-invertible, with $M\geq 3$. Thus $(X_{\btau_{1}},\sigma)$ is
quasi-invertible, even though $(X_{\btau},\sigma)$ is not, and we
cannot generate an adic representation for $\btau$ from one for
$\btau_{1}$. For example, if $\btau(0)=00100110$, then the return
words are ${\cal R}=\{0,01,011\}$ and $\btau_{1}$ is 3-almost
periodic, with $\bx^{1}, \bx^{2}, $ and $\bx^{3}$ the preimages of
${\cal D} (\bu)$. In $X_{\btau}$ though, $\phi(\bx^{1})$, and $\sigma
(\phi(\bx^{2}))$ are preimages of $\bu$, and $\phi(\bx^{2})$ and
$\sigma(\phi(\bx^{3}))$ are preimages of $\sigma (\phi(\bx^{2}))$. So
an adic representation of $(X_{\btau},\sigma)$ would need 4 maximal
elements, one of which is also minimal, and another minimal element.

\item\label{orbitequivalence}
 For any primitive substitution $\btau$ whose composition matrix has a
 rational Perron-Frobenius eigenvalue, there exists a quasi-invertible
 primitive substitution $\btau'$ whose associated two sided substitution system 
$(X_{\btau}^{\mathbb Z},\sigma)$
 is orbit equivalent to $(X_{\btau}^{\mathbb Z}, \sigma)$. To see this we use the
 results in Corollary 6.7 and Theorem 6.15 in \cite{y}. First, given
 $\btau$, with Perron value $\lambda$, we find a positive integer $d$
 so that $< \mu (E): E\subset X_{\btau}\mbox{ is clopen } >= \{ \frac{n}{m}: n\,\in
 {\mathbb Z}, m \mbox{ is a factor of } d.\lambda^{n} \mbox {for some }
 n\}$, where the former is order isomorphic to the dimension group modulo the infinitesimal subgroup of $(X_{\btau}, \sigma)$. Then given $d$ and $\lambda$, both greater than 1, the author
 describes conditions on the composition matrix of ${\btau'}$ so that
 $< \mu (E): E\subset X_{\btau'}, E \mbox{ clopen }>= \{ \frac{n}{m}: n\,\in {\mathbb
 Z}, m \mbox{ is a factor of } d.\lambda^{n} \mbox { for some } n\}$.  In particular if the substitution $\btau'$ defined on $d$ letters has the substitution matrix 

\[ \left( \begin{array}{ccccc}
\lambda^{m}& \lambda^{m} & \ldots & \lambda^{m} &(\lambda^{m}-d+1)\lambda^{m} \\
\lambda^{m}&\lambda^{m}  & \ldots & (\lambda^{m}-d+1)\lambda^{m} & \lambda^{m}\\\vdots & \vdots & \vdots &\vdots &\vdots\\
(\lambda^{m}-d+1)\lambda^{m}&\lambda^{m}  &\ldots &\lambda^{m} &\lambda^{m} \end{array} \right)\] 
then $\btau'$ has the desired dimension group (modulo the infinitesimal subgroup)
 We
 claim here that $\btau'$ can be chosen to be quasi invertible. In particular if $\btau $ is defined by 
$\btau'(1)=1^{\lambda^{m}}2^{\lambda^{m}}3^{\lambda^{m}} \ldots (d-1)^{\lambda^{m}}
d^{({\lambda^{m}}-d)\lambda^{m}} d^{{\lambda^{m}}}$,
$\btau'(2)=1^{\lambda^{m}}2^{\lambda^{m}}3^{\lambda^{m}} \ldots (d-1)^{\lambda^{m}}
(d-1)^{({\lambda^{m}}-d)\lambda^{m}} d^{{\lambda^{m}}}, \ldots 
\btau'(d)=1^{\lambda^{m}}2^{\lambda^{m}}3^{\lambda^{m}} \ldots (d-1)^{\lambda^{m}}
(1)^{({\lambda^{m}}-d)\lambda^{m}} d^{{\lambda^{m}}},$
then $\btau'$ is quasi invertible. 

\item\label{notonesidedconjugate}
 Here is an example of a quasi-invertible substitution whose two-sided
subshift conjugacy class is not the same as its one-sided subshift
conjugacy class.  Define $\btau (a) = aabaa$, $\btau (b) = abcab$, and
$\btau (c) = aabac$.  The fixed point $\bu = a\ldots$ is a 3-branch
point, and $(X_{\btau}^{\mathbb N}, \sigma)$ is quasi-invertible. Let
$\bu_{L}=\ldots u_{-2}u_{-1}$ be the left infinite sequence generated
by $\{\btau^{n}(a)\}$, ${\bf v}_{L}$ the one generated by
$\{\btau^{n}(b)\}$, and ${\bf w}_{L}$ the one generated by
$\{\btau^{n}(c)\}$. Note that if $n\neq -1$, $u_{n}=w_{n}$. The words
of length 3 in ${\cal L}_{\btau}$ are ${\cal
W}_{3}=\{aab,\,\,\,aba,\,\,\,baa,\,\,\,aaa,\,\,\,abc,\,\,\, bca,\,\,\,
cab,\,\,\,bac,\,\,\,aca,\,\,\,caa\}$. Define a local rule $\phi:{\cal
W}_{3}\rightarrow \{\alpha, \beta, \gamma, d,e,f,g,h\}$ with left and
right radius one, as $\phi (aaa)=\phi(baa)=\phi (aca) = \alpha$,
$\phi(aba)=\beta$, $\phi(caa)=\gamma$, and let $\phi$ map the
remaining words in ${\cal W}_{3}$ in a one-one fashion to
$\{d,\,\,e,\,\,f,\,\,g,\,\,h\}$. Let $\Phi$ be the shift-commuting
factor mapping corresponding to $\phi$ and Let $Y=\Phi(X_{\btau}^{\mathbb Z})$.

\begin{lemma}
The map $\Phi:X_{\btau}^{\mathbb Z}\rightarrow Y$ is injective.
\end{lemma}
{\bf Proof:} Assume that for some ${\bf x}\neq {\bf y}$, $\Phi({\bf
x})= \Phi({\bf y})$. So there is some $n$ such that
$x_{n-1}x_{n}x_{n+1}\neq y_{n-1}y_{n}y_{n+1}$, and yet both words
belong to $\{baa,\,\, aaa,\,\,aca\}$. There are three cases to
consider, all similar to prove, so we look at the case when
$x_{n-1}x_{n}x_{n+1}=baa$ and $y_{n-1}y_{n}y_{n+1}=aaa$. If $baa$
appears in a sequence here, it either appears as a subword of
$\btau(a)= aabaa$, or $\btau(ba)=abcabaabaa$ or
$\btau(bc)=abcabaabac$. On the other hand $aaa$ can only appear as a
subword of $\btau(aa)=aabaaaabaa$ or $\btau(ac)=aabaaaabac$ or 
$\btau(ab)=aabaaabcab$. This
means that $x_{n-2}x_{n-1}x_{n} = aba$ and $y_{n-2}y_{n-1}y_{n}=aaa$ or $baa$. Either way 
$\Phi({\bf x})_{n}= \beta$ and $\Phi({\bf y})_{n}=\alpha$, contradicting $\Phi({\bf
x})= \Phi({\bf y})$.
\endproof

Now Theorem 4 in \cite{dhs} tells us that $(Y,\sigma)$ is either a
substitution system,or a stationary odometer. Since $(X_{\btau},\sigma)$ is not an odometer, and $\Phi$ is a conjugacy, so $(Y,\sigma)$  is generated by some substitution
$\btau'$. Hence the 2-sided subshifts generated by $\btau$ and $\btau'$ are conjugate. 

If ${\bf x}=x_{0}x_{1}\ldots$ is a one sided sequence, let $\Phi(x)= \phi(x_{0}x_{1}x_{2}) \phi(x_{1}x_{2}x_{3})\ldots$.
Note that $\btau'$ has two 2-branch points. For, \begin{eqnarray*}\Phi({\bf
u}_{L}\cdot {\bf u} ) = \Phi(\ldots u_{-2}u_{-1})\phi(u_{-2}u_{-1}u_{0})\cdot
\phi(u_{-1}u_{0}u_{1})\Phi ({\bf u}) =\\ \Phi(\ldots u_{-2}u_{-1})\phi(aaa)\cdot
\phi(aaa)\Phi ({\bf u}) = \Phi(\ldots u_{-2}u_{-1})\alpha\cdot
\alpha\Phi ({\bf u})\end{eqnarray*}
and 
 \begin{eqnarray*}\Phi({\bf
w}_{L}\cdot {\bf u} ) = \Phi(\ldots w_{-2}w_{-1})\phi(w_{-2}w_{-1}u_{0})\cdot
\phi(w_{-1}u_{0}u_{1})\Phi ({\bf u}) = \\\Phi(\ldots w_{-2}w_{-1})\phi(aca)\cdot
\phi(caa)\Phi ({\bf u}) = \Phi(\ldots w_{-2}w_{-1})\alpha\cdot
\gamma\Phi ({\bf u}),\end{eqnarray*}
so that $\Phi ({\bf u})$ is a 2-branch point; while 
\begin{eqnarray*}\Phi({\bf
v}_{L}\cdot {\bf u} ) = \Phi(\ldots v_{-2}v_{-1})\phi(v_{-2}v_{-1}u_{0})\cdot
\phi(v_{-1}u_{0}u_{1})\Phi ({\bf u}) = \\\Phi(\ldots v_{-2}v_{-1})\phi(aba)\cdot
\phi(baa)\Phi ({\bf u}) = \Phi(\ldots v_{-2}v_{-1})\beta\cdot
\alpha\Phi ({\bf u}),\end{eqnarray*}
so that $\alpha\Phi({\bf u})$ is also a 2-branch point. Thus the one sided substitution systems generated by $\btau$ and $\btau'$ cannot be topologically conjugate.

\end{enumerate}

\section{(N,n)-substitutions.}

In this section we consider substitutions which have more than one branch point. The next theorem is a special case of Theorem \ref{generalfixedbranch}, which we include as its proof is simpler:

\begin{theorem}
\label{specialfixedbranch}
Suppose that $\btau$ is defined on ${\cal A}=
\{0,1,\ldots n-1\}$, and is prefix and suffix permutative, so that
there are $n$ fixed points, the $i$-th of which is a $k_{i}$-branch
point. If ${\cal L}_{\btau} $ has $N$ words of length 2, then there
exists $\btau'$ defined on ${\cal A}'=\{0,1,\ldots N-1\}$ such that
$\btau'$ is suffix permutative, with $n$ fixed points, where each is a
$k_{n}$ branch point and where $(X_{\btau},\sigma)\cong(X_{\btau'},\sigma)$, which is topologically conjugate to 
$(X_{{\cal B}_{\btau'}},V_{\btau'})$. 
\end{theorem}

{\bf Proof:} First note  that all branch points are $\btau$-fixed. By taking a power of $\btau$ if necessary we can assume that $\btau(i)$ begins and ends with $i$.
 Ordering the words of length two in ${\cal
L}_{\btau}$ lexicographically, define the right radius one local
rule $\phi$ on these words by mapping the $j$-th element in this list to $j$, and let $\Phi$ be the shift commuting map associated with $\phi$.
Note that the image of $\phi$ is ${\cal A}':=\{0,\ldots ,N-1\}$ where 
 $N$ is  the number of words of length 2 in ${\cal L}_{\btau}$ and satisfies $N=\sum_{i=0}^{n-1}k_{i}$.
 Now
define the substitution $\btau'$ on ${\cal A'}$ as follows:
if $a=\phi(ij)$, then ${\btau'}(a):= \phi (\btau(i))a$. It is clear
that $\btau'$ is suffix permutative, with $n$ fixed points. We will
show that $(X_{\btau},\sigma)$ is topologically conjugate
to $(X_{\btau'} ,\sigma)$ by showing that
${\bf v_{i}}:=\Phi(\btau^{\infty}(i))$ is the $\btau'$-fixed in
$(X_{\btau'} ,\sigma)$. If that is the case, then since
$\Phi$ is injective and  shift commuting, and the fixed points
generate the subshift, we are done.

To see this we show by induction that $\btau'^{n}(a)=
\phi(\btau^{n}(i))a$ for all $n$, whenever $\phi(ij)=a$. For $n=1$ the
assertion is true by definition. Assuming the case true for $n$, we have
\begin{eqnarray*}
& & \btau'^{n+1}(a) 
=  \btau'^{n}(\btau'(a)) = \btau'^{n} (\phi(\btau(i)) a)
 =  \btau'^{n}(\phi(ii_{1}\ldots i_{l_i}i)a)\\ &  = &
\btau'^{n} \left(\phi(ii_{1})\phi(i_{1}i_{2})\ldots \phi(i_{l_{i}}i)a\right)\\
& = & 
\btau'^{n}(\phi(ii_{1})) \btau'^{n}(\phi(i_{1}i_{2}))\ldots \btau'^{n}(\phi(i_{l_{i}}i)) \btau'^{n}(a)\\
 & = &
\phi(\btau^{n}(i)) \phi(ii_{1})
\phi(\btau^{n}(i_{1})) \phi(i_{1}i_{2}) \ldots 
 \phi(\btau^{n}(i_{l_i}))  
\phi(i_{l_{i}}i)
\phi({\btau}^{n}(i))a \\
& \stackrel{*}{=} &
\phi\left(\btau^{n}(i)
\btau^{n}(i_{1})\ldots \btau^{n}(i_{l_{i}}) \btau^{n}(i)           
\right)a\, 
\\ & = & \phi\left(\btau^{n+1}(i)\right)a\, ,\\
\end{eqnarray*}
 
where * follows since $\btau$ is prefix and suffix permutative.
Now apply this to  ${\bf v_{i}} \,\in \phi[ni+i_{1}]$.

To see that $\btau'$ has an adic representation, note that if
$i'=\phi(ij)$ generates a $\btau'$-fixed point, then the $k_{i}$
predecessors of $i'$ are $\{\phi(\ell i)\}$ as $\ell$ ranges over all
predecessors of $i$ in ${\cal L}_{\btau}$. As for distinct $i$, these
sets are disjoint, Proposition \ref{closingproposition} implies that
$V_{\btau'}$ is continuous, and Proposition \ref{rightadicrep} that $(X_{\btau'} ,\sigma)\cong
(X_{{\cal B}_{\btau'}},V_{\btau'})$.

\endproof
{\bf Examples:} 
\begin{enumerate}
\item
If $\btau (0)=01$ and $\btau(1)=10$, then $\btau'$ is
defined as $\btau'(0) = 1320$, $\btau'(1)=1321$, $\btau' (2) = 2012$,
and $\btau' (3) = 2013$. Note that ${\cal B}_{\btau'}$ has four
maximal elements $M_{0},\ldots ,M_{3}$ and two minimal elements
$m_{1},m_{2}$, where $M_{i}$ is constant in ${\cal B}_{\btau}$
through vertex $i$.
The words of length 2 in ${\cal L}_{\btau'}$ are $
\{ 01,\,12,\, 13,\,20,\,21,\,32\}$ so $V_{\cal B}:X_{\cal
B}\rightarrow X_{\cal B}$ can be defined continuously with $V_{\cal
B}(M_{0}) = V_{\cal B}(M_{2})= m_{1}$ and $V_{\cal B}(M_{1})=V_{\cal
B}(M_{3}) = m_{2}$.
\item
 Let $\btau(0)=0120$, $\btau(1) = 100101121$, $\btau (2) =
 20212$. Then $\btau^{\infty}(0)$ and $\btau^{\infty}(1)$ are 3-branch
 points, and $\btau^{\infty}(2)$ is a 2-branch point. Here $N=8$, so
 define $\btau'$ on $\{0,1,2,\ldots 7\}$ by $\btau' (0)=1560$,
 $\btau'(1)=1561$, $\btau' (2)= 1562$, $\btau' (3) = 301314573$,
 $\btau'(4) = 301314574$, $\btau'(5) = 3013141575$, $\btau'(6) =
 62756$, and $\btau'(7)=62757$.  One can check that
 ${\btau'}^{\infty}(1)$ is a 3-branch point with predecessors
 $0{\btau'}^{\infty}(1)$, $3{\btau'}^{\infty}(1)$, and
 $6{\btau'}^{\infty}(1)$; that ${\btau'}^{\infty}(3)$ is also a
 3-branch point with predecessors $1{\btau'}^{\infty}(3)$,
 $4{\btau'}^{\infty}(3)$ and $7{\btau'}^{\infty}(3)$, and that
 ${\btau'}^{\infty}(6)$ is a 2-branch point with predecessors
 $2{\btau'}^{\infty}(6)$ and $7{\btau'}^{\infty}(6)$. Thus extend
 $V_{{\cal B}_{\btau'}}$ by mapping $V_{{\cal B}_{\btau'}}({\bf
 M_{0}})=V_{{\cal B}_{\btau'}}({\bf M_{3}})= V_{{\cal
 B}_{\btau'}}({\bf M_{6}})= {\bf m_{0}}$, $V_{{\cal B}_{\btau'}}({\bf
 M_{1}})=V_{{\cal B}_{\btau'}}({\bf M_{4}})= V_{{\cal
 B}_{\btau'}}({\bf M_{7}})= {\bf m_{1}}$ and $V_{{\cal
 B}_{\btau'}}({\bf M_{2}})=V_{{\cal B}_{\btau'}}({\bf M_{5}})= {\bf
 m_{2}}$.
\end{enumerate}

\begin{theorem}\label{generalfixedbranch}
Suppose that $\btau$ is defined on ${\cal A}= \{0,1,\ldots n-1\}$, so
that there are $m$ right-fixed points $\{{\bf p}_{i}: 1\leq i \leq m \}$,
each of which is a $k_{i}$-branch point. Suppose there are no other branch points. If ${\cal L}_{\btau} $ has
$N$ words of length 2, then there exists $\btau'$ defined on the
$N$-letter alphabet ${\cal A}'$ such that 
\begin{enumerate}
\item
$(X_{\btau},\sigma) \cong (X_{\btau'},\sigma)$,
\item
$\btau'$ has  $m$ right-fixed points $\{{\bf m}_{i} = (m_{i}, \ldots ):
1\leq i \leq m\},$ 
\item $\btau'$ has  $\sum_{i=1}^{m}k_{i}$ left-fixed points, 
\item $V_{\btau'}$ is continuous, and 
\item $(X_{\btau'} ,\sigma)\cong (X_{{\cal B}_{\btau'}},V_{\btau'})$ 
\end{enumerate}
\end{theorem}

    {\bf Proof:} Assume that the fixed
point ${\bf p_{i}}$ starts with the letters $p_{i}p_{i}^{1}$. We can assume, by taking a power of $\btau$ if
necessary, that all substitutions words start with $p_{i}p_{i}^{1}$ for some 
$i\,\in \{1,\ldots m\}$, and end with a letter
which generates a left infinite $\btau$-fixed point. We define a
right-sided, radius one factor mapping $\phi:{\cal L}_{2}(\btau)
\rightarrow {\cal A}'$ where ${\cal A}'$ is an $N$-letter alphabet. 
Let $\phi(p_{i},p_{i}^{1}):= m_{i}$. If
$s{\bf p_{i}}$ is a $\sigma$-predecessor of ${\bf p}_{i}$, where the
cylinder set $[sp_{i}]$ has not appeared as the home of one of the
right infinite $\btau$-fixed points, define $\phi(sp_{i}) = M_{l}$ for some
$l$, as $p_{i}$ and $s$ range over all allowed letters. For the remaining
words $ab$ in ${\cal L}_{2}(\btau)$, let $\phi$ map them injectively
to letters that are distinct from the already specified $m_{i}$'s and
$M_{i}$'s. Let ${\cal A}'$ be the resulting alphabet.  Define $\tau'$
on ${\cal A}'$ as follows. If $i\,\in {\cal A}$, let $p^{i}$ be the first
letter of $\btau(i)$, and $s^{i}$ the last letter of $\btau(i)$. If $a = \phi(ij)$, let $\btau'(a) :=
\phi(\btau(i))\phi(s^{i}p^{j})$. Here $\phi(\btau(i))$ is a word of
length $|\btau(i)|-1$.   The
following two claims are straightforward to verify, starting with the
definition of $\btau'$.

{\bf Claim 1:} If $\alpha \beta$ is a 2-word appearing in
$\btau'(\gamma)$, where $\alpha=\phi(ij)$ and $\beta = \phi(i'j')$,
then the suffix of $\btau'(\alpha)$ is $s=\phi(s^{i}p^{i'})$, as is the suffix of $\btau'^{n}(\alpha)$ for each $n\,\in {\mathbb N}$.


{\bf Claim 2:} If $\alpha = \phi (i_{\alpha}, j_{\alpha})$ for 
$\alpha\, \in{\cal A}'$, and $\btau'(\alpha)= 
\alpha_{1}\,\alpha_{2}\, \ldots \,\alpha_{k}$, then $\btau(i_{\alpha}) = 
i_{\alpha_{1}}i_{\alpha_{2}}\ldots i_{\alpha_{k}}$.


{\bf Claim 3:} If $\alpha = \phi(i_{\alpha},j_{\alpha})$ is a letter in ${\cal A'}$, then
$(\btau')^{n}(\alpha)\backslash \mbox{ a suffix }= \phi
(\btau^{n}(i))$ for each $n\, \in {\mathbb N}$. If $n=1$ then it is true by definition.

If it is true for $n$, then 
\begin{eqnarray*}
(\btau')^{n+1}(\alpha))& = &
(\btau')^{n}(\btau'(\alpha)) = (\btau')^{n}(\alpha_{1}\alpha_{2}\ldots\alpha_{k})
\\& = &(\btau')^{n}(\alpha_{1})\, (\btau')^{n}(\alpha_{2})\ldots\,(\btau')^{n}(\alpha_{k})\\ &\stackrel{\mbox{Claim }1}{=}&
\phi(\btau^{n}(i_{\alpha_{1}}))\,\phi(s^{i_{\alpha_{1}}}, p^{i_{\alpha_{2}}})\,
\phi(\btau^{n}(i_{\alpha_{2}}))\,\phi(s^{i_{\alpha_{2}}}, p^{i_{\alpha_{3}}})\,\ldots
\phi(\btau^{n}(i_{\alpha_{k}}))\ell
\\
&=&
\phi(\btau^{n}(i_{\alpha_{1}}i_{\alpha_{2}}\ldots i_{\alpha_{k}}))\ell {\stackrel{\mbox{Claim }2}{=}} 
\phi(\btau^{n}(\btau(i_{\alpha}))\ell = \phi(\btau^{n+1}(i_{\alpha}))\ell\, ,
\end{eqnarray*}
where $\ell$ is some letter in ${\cal A}'$. The result follows.

Thus in particular $\Phi$ maps the $\btau$-fixed points $\{{\bf
p}_{i}\} $ to the $\btau'$-fixed points $\{{\bf m}_{i}\}$ and since
$\phi$ is 1-1, $\Phi$ is shift commuting, and the fixed points
generate the subshift, then $(X_{\btau},\sigma)$ is
topologically conjugate to $(X_{\btau'}^{\mathbb
N},\sigma)$. Statements 2 and 3 follow by definition of $\btau'$, and to
see 4, note that if $m_{k}\neq m_{l}$, then
$i_{m_k}\neq i_{m_l}$.Thus if $Mm_{k}\,\in {\cal L}_{\btau'}$, then
$Mm_{i}\,\not \in {\cal L}$ for any other $m_{i}$. Thus $M$ can be
followed by at most one of the $m_{i}'s$, and must be followed by at
least one, by minimality of $\btau$.

Clearly $\btau'$ is primitive. We show that $\btau'$ is recognizable and
$\btau':X_{\btau'} \rightarrow X_{\btau'} $ is a homeomorphism (here
$\btau'$ is not injective, so that Part 4 of Proposition
\ref{hosttopologyprop} cannot be directly applied). If  $\bx'$, $\by'\, \in X_{\btau'} $,
let the $\Phi$-preimages of $\bx'$ and $\by'$ in $X_{\btau}$ be $\bx$
and $\by$ respectively.  If $\btau'(\bx') = \btau'(\by')$, then using Part
3 of Proposition \ref{hosttopologyprop}, we have $|\btau'(x_{n}')| =
|\btau'(y_{n}')|$ for each $n\geq 0$, so that $\btau'(x_{n}') =
\btau'(y_{n}')$ for each $n\geq 0$. Now $\btau' (x_{n}') =
\phi(\btau(x_{n}))\phi(s_{x_n}p_{x_{n+1}})$ and $\btau' (y_{n}') =
\phi(\btau(y_{n}))\phi(s_{y_n}p_{y_{n+1}})$, and since $\btau$ is
injective, if $\bx_{n}\neq \by_{n}$, either $\phi(\btau(x_{n}))\neq
\phi(\btau(y_{n}))$, or $s_{x_n}\neq s_{y_n}$. Thus $\bx=\by$ and so
$\bx'=\by'$. If $\btau'$ is not recognizable, then there exist letters
$a'$, $b'$ and for each $L$, words $w_{L}'$ satisfying the conditions
of Theorem \ref{mosserecognizability}. This implies that there exist
corresponding letters $a$ $b$ such that $\btau(b)$ is a proper suffix of $\btau(a)$, and, as in the proof of injectivity of $\btau'$  the words $w_{L}'$ lift to words 
$w_{L}\,\in {\cal L}_{\btau}$ such that conditions of  Theorem \ref{mosserecognizability}
are satisfied for $\btau$ contradicting recognizability of $\btau$.

Thus all branch points of $\btau'$ are fixed points, and using  Corollary \ref{clopenpartition}, the partitions $({\cal P}_{n}')$ generate the topology of 
$X_{\btau'} $. Now  Proposition \ref{rightadicrep} implies that  $(X_{\btau'} ,\sigma)\cong (X_{{\cal B}_{\btau'}},V_{\btau'})$. 
\endproof

{\bf Example:} If $\btau(0)=01200, \btau(1)=11201$ and $\btau(2) =
10121$, then $\btau$ has two fixed points, both of which are 2-branch
points. Define $\phi(01)=m_{1}$, $\phi(11)=m_{2}$, $\phi(00)=M_{1}$,
$\phi(10)=M_{2}$, $\phi(12)=a$, $\phi(20)=b$ and $\phi(21)=c$. If
$A_{i}:=\phi(\btau(i))$, then $A_{0}= m_{1}abM_{1}$, $A_{1}=
m_{2}abm_{1}$ and $A_{2} = M_{2}m_{1}ac$; and $\btau'$ is defined as
$\btau'(M_{1})=A_{0}M_{1}$, $\btau'(m_{1}) = A_{0}m_{1}$,
$\btau'(M_{2})= A_{1}M_{2}$, $\btau(m_{2})=\btau(a)= A_{1}m_{2}$,
$\btau(b)=A_{2}M_{2}$ and $\btau(c)=A_{2}m_{2}$.  Note $\btau'$ has
four left infinite fixed points, generated by $\{m_{1}, m_{2}, M_{1},
M_{2}\}$, and the predecessors of $m_{1}$ are $\{M_{1},M_{2}\}$, while
those of $m_{2}$ are $\{m_{1}, m_{2}\}$.

{\footnotesize \bibliographystyle{alpha} \bibliography{bibliography} }

\end{document}